\documentclass[12pt,a4paper,reqno]{amsart}
\usepackage{amssymb}
\usepackage{multirow}
\usepackage{graphicx}

\usepackage[main=english]{babel}
\newtheorem{Definition}{Definition}[section]
\newtheorem{Example}{Example}[section]
\newtheorem{Theorem}{Theorem}[section]
\newtheorem{Remark}{Remark}[section]
\newcommand{\D}{\mathbb{D}}

\newcommand{\I}{\mathbb{I}}
\newcommand{\R}{\mathbb{R}}

\newcommand{\Com}{\mathbb{C}}
\newcommand{\N}{\mathbb{N}}

\begin{document}


\title[Fractional Differential Equations]{Fractional Differential Equations with the General Fractional Derivatives of Arbitrary Order in the Riemann-Liouville Sense}

\author{Yuri Luchko}
\curraddr{Department of Mathematics, Physics, and Chemistry, Berlin University of Applied Sciences and Technology,  
     Luxemburger Str. 10,  
     13353 Berlin,
Germany}
\email{luchko@bht-berlin.de}

\subjclass[2010]{26A33;  26B30; 33E30; 44A10; 44A35; 44A40;  45D05; 45E10; 45J05}
\dedicatory{}
\keywords{Sonine kernel; Sonine condition;  general fractional integral; general fractional derivative of arbitrary order; fundamental theorems of fractional calculus; operational calculus; fractional differential equations; convolution series}

\begin{abstract}
In this paper, we first deal with the general fractional derivatives of arbitrary order defined in the Riemann-Liouville sense. In particular, we deduce an explicit form of their null space and prove the second fundamental theorem of Fractional Calculus that leads to a closed form formula for their projector operator. These results allow us to formulate the natural initial conditions for the fractional differential equations with the 
general fractional derivatives of arbitrary order in the Riemann-Liouville sense. In the second part of the paper, we develop an operational calculus of the  Mikusi\'nski type for the general fractional derivatives of arbitrary order in the Riemann-Liouville sense and apply it for derivation of an explicit form of solutions to the Cauchy problems for the single- and multi-term linear fractional differential equations  with these derivatives. The solutions are provided in form of the convolution series generated by the kernels of the corresponding general fractional integrals.
\end{abstract}

\maketitle


\section{Introduction}
\label{sec1}

In the framework of the abstract Volterra integral equations on the Banach spaces, the evolution equations with the integro-differential operators of the convolution type
with different classes of kernels and on different spaces of functions were a subject of active research within several last decades (see \cite{Cle84}, \cite{Pr} and references therein). This research can be traced back to the paper \cite{Son84} by Sonine published in 1884, who extended the Abel's method for derivation of an analytical solution to the Abel integral equation to a quite general class of the integral equations in the form 
\begin{equation}
\label{Son1}
f(t) = (\kappa\, *\, \phi)(t)  \, = \, \int_0^t \kappa(t-\tau)\phi(\tau)\, d\tau,\ t>0
\end{equation}
with the kernels $\kappa$  that satisfy the Sonine condition formulated as follows: A function $\kappa$ is called a Sonine kernel if there exists a function $k$ such that 
\begin{equation}
\label{Son}
(\kappa \, *\, k )(t) = \int_0^t \, \kappa(t-\tau)\, k(\tau)\, d\tau\ = \ \{1 \},\ t>0,
\end{equation}
where by $*$ we denote the Laplace convolution and by $\{1 \}$ the function identically equal to one for $t>0$. For a Sonine kernel $\kappa$, the kernel $k$ from the relation \eqref{Son} is called its associated kernel. In \cite{Sam}, several examples of the Sonine kernels and their properties can be found.

For this class of kernels,  Sonine showed that the integral equation \eqref{Son1}
possesses a (formal) solution in the form
\begin{equation}
\label{Son2}
\phi(t) = \frac{d}{dt} (k\, *\, f)(t)  \, = \, \frac{d}{dt}\, \int_0^t k(t-\tau)f(\tau)\, d\tau,\ t>0.
\end{equation}

In the case of the Abel integral equation  
\begin{equation}
\label{Abel1}
f(t) = \frac{1}{\Gamma(\alpha)}\, \int_0^t (t-\tau)^{\alpha-1}\phi(\tau)\, d\tau, \ t>0,\ 0<\alpha <1
\end{equation}
that was treated by Abel in the papers \cite{Abel1,Abel2}, the kernel $\kappa(t) = \frac{t^{\alpha-1}}{\Gamma(\alpha)}$ is a Sonine kernel with an associated kernel $k(t) = \frac{t^{-\alpha}}{\Gamma(1-\alpha)}$ because of the well-known relation
\begin{equation}
\label{Abel2}
(h_\alpha \, *\, h_{1-\alpha} )(t) = \ \{1 \},\ t>0,\ 0<\alpha <1,
\end{equation}
where by $h_\alpha$ we denote the following power function:
\begin{equation}
\label{h}
h_\alpha(t) = \frac{t^{\alpha-1}}{\Gamma(\alpha)},\ t>0,\ \alpha >0.
\end{equation}
Thus, the solution to the Abel integral equation \eqref{Abel1} takes the well-known form:
\begin{equation}
\label{Abel3}
\phi(t) = \frac{d}{dt}\, \frac{1}{\Gamma(1-\alpha)}\, \int_0^t (t-\tau)^{-\alpha}f(\tau)\, d\tau,\ t>0.
\end{equation}

In the modern Fractional Calculus (FC), the right-hand side of the Abel integral equation \eqref{Abel1} is referred to as the Riemann-Liouville fractional integral of the order $\alpha,\ \alpha >0$:
\begin{equation}
\label{RLI}
(I^\alpha_{0+}\, f)(t) := (h_\alpha\, * \, f)(t) \, = \, \frac{1}{\Gamma(\alpha)}\, \int_0^t (t-\tau)^{\alpha-1}f(\tau)\, d\tau,\ t>0,
\end{equation}
whereas the right-hand side of the Abel solution formula \eqref{Abel3} is called the Riemann-Liouville fractional derivative of the order $\alpha,\ 0\le \alpha<1$:
\begin{equation}
\label{RLD1}
(D^\alpha_{0+}\, f)(t) := \frac{d}{dt}(h_{1-\alpha}\, * \, f)(t) \, = \,\frac{d}{dt}\, (I^{1-\alpha}_{0+}\, f)(t),\ t>0.
\end{equation}

Because of the evident connection of the Riemann-Liouville fractional integral and derivative to the operators \eqref{Son1} and \eqref{Son2} with the Sonine kernels, an interpretation of these general operators in the framework of FC was just a question of time. The first publication entirely devoted to this interpretation was the paper \cite{Koch11}, where the operator
\begin{equation}
\label{FDR-L}
(\D_{(k)}\, f)(t) = \frac{d}{dt} (k\, *\, f)(t) = \frac{d}{dt}\, \int_0^t k(t-\tau)f(\tau)\, d\tau,\ t>0
\end{equation}
with the Sonine kernels from a special class  $\mathcal{K}$ was called the general fractional derivative (GFD) of the Riemann-Liouville type.  The GFD of the Caputo type was defined in \cite{Koch11} in the form
\begin{equation}
\label{FDC}
( _*\D_{(k)}\, f) (t) =  (\D_{(k)}\, f) (t) - f(0)k(t),\ t>0.
\end{equation}
In \cite{Koch11}, mainly the GFD of the Caputo type was treated. In particular, the GFD \eqref{FDC} with the kernels from  $\mathcal{K}$ was shown to be a left inverse operator to the general fractional integral (GFI)
\begin{equation}
\label{GFI}
(\I_{(\kappa)}\, f)(t) = (\kappa\, *\, f)(t) \, = \, \int_0^t \kappa(t-\tau)f(\tau)\, d\tau,\ t>0,
\end{equation}
where $\kappa$ is the Sonine kernel associated to the kernel $k$ of the GFD \eqref{FDC}.  Moreover, in \cite{Koch11}, some important properties of  the solutions to the fractional relaxation equation and to  the Cauchy problem for
the  time-fractional diffusion equation with the GFD \eqref{FDC} were investigated in detail. 

In the next paper \cite{LucYam16} devoted to the GFD \eqref{FDC}, a maximum principle for the initial-boundary-value problems for the time-fractional diffusion equations with the GFD of Caputo type  was deduced. Further results regarding  the ordinary and partial time-fractional  differential equations with the GFD \eqref{FDC} have been derived in \cite{KK,Koch19_1,Koch19_2,LucYam20,Sin18,Sin20}. For results regarding inverse problems for the fractional differential equations involving the GFD and their applications we refer to \cite{JK17,KJ19_1,KJ19}.

Very recently, a series of papers devoted to the GFI and GFD with the Sonine kernels that possess an integrable singularity of power function type at the point zero has been published. In \cite{Luc21a}, these operators were studied on the space of functions continuous on the real positive semi-axis that can have an integrable singularity of power function type at the point zero. In \cite{Luc21b}, the  GFDs of arbitrary order in the Riemann-Liouville and in the Caputo senses were introduced and investigated. These operators extend the definitions \eqref{FDR-L} and \eqref{FDC} that correspond to the case of the "generalized order" from the interval $(0,\, 1)$ to the case of any positive real order. An important subclass of the kernels of the GFDs of arbitrary order was suggested in \cite{Tar}. Operational method for solving the Cauchy problems for the fractional differential equations with the GFDs of the Caputo type was worked out in \cite{Luc21c}. In \cite{Luc21d}, analytical solutions to some fractional differential equations with the GFDs of the Riemann-Liouville type with the suitable formulated initial conditions were derived for the first time in the FC literature. The solution method employed in \cite{Luc21d} is the technique of the convolution series that are a far reaching generalization of the power series with both integer and fractional exponents. Finally, we mention the papers \cite{Tar1,Tar2,Tar3}, where the theory presented in \cite{Luc21a,Luc21b,Tar,Luc21c,Luc21d} was applied for formulation of a general  fractional dynamics, a general non-Markovian quantum dynamics, and a general fractional vector calculus.

For applications, the fractional differential equations with different kinds of the fractional derivatives are especially important. As in the case of the differential equations with the integer order derivatives, the class of the fractional differential equations that can be solved in explicit form is very restricted. These are mainly  linear equations with the constant coefficients. However, exactly this class of equations is especially important both for applications and as a basis for a qualitative treatment of quasi-linear and non-linear equations. One of the most simple and powerful techniques for analytical treatment of the linear ordinary and fractional differential equations is by means of the operational calculi of Mikusi\'nski type developed for the corresponding derivatives.   For a presentation of an operational calculus for the first order derivative and its applications we refer to \cite{Mik59}.   An operational calculus of the Mikusi\'nski type for the hyper-Bessel differential operator has been suggested in \cite{Dim66}. 

In \cite{Luc93}, the first operational calculus for a fractional derivative has been developed and applied for analytical treatment of the fractional differential equations. This operational calculus was constructed for the multiple  Erd\'elyi-Kober fractional derivative (see also \cite{LucYak94,YakLuc94,BasLuc95}) and applied for analytical treatment of the fractional differential equations involving this derivative. The case of the Riemann-Liouville fractional derivative was treated in \cite{LucSri95,HadLuc}. Operational method for derivation of the closed form solutions to the initial-value problems for the single- and multi-term fractional differential equations was presented in \cite{LucSri95} for the case of commensurate and in \cite{HadLuc} for the case of non-commensurate orders of the Riemann-Liouville fractional derivatives. In \cite{LucGor99}, an operational calculus for the Caputo fractional derivative was developed and applied for 
solving fractional differential equations involving the Caputo fractional derivatives with the commensurate and non-commensurate orders. Other operational calculi of Mikusi\'nski type were developed  in \cite{HLT09} for the generalized Riemann-Liouville fractional derivative (Hilfer fractional derivative), in \cite{Han}
 for the Caputo-type fractional Erd\'elyi-Kober  derivative,  and in \cite{Fer1,Fer2} for the Riemann-Liouville and Caputo fractional derivatives with respect to functions, respectively. In \cite{GorLuc97}, the 
operational method was applied for solving the 
generalized Abel integral equations of the second kind and in \cite{GLS97}, the integral equations with the Gauss hypergeometric function as a kernel have been solved by means of a suitably constructed operational calculus. In the recent paper \cite{Luc21c}, a  Mikusi\'nski type operational calculus was developed for the GFD in the Caputo sense  with the Sonine kernels that possess an integrable singularity of power function type at the point zero. As an application, this operational calculus was employed to derive analytical solutions to the Cauchy problems for the single- and multi-term fractional differential equations with the GFDs in the Caputo sense. Finally, we mention the surveys  \cite{Luc99,Luc19}, where several operational calculi for different fractional derivatives as well as operational method for solving the fractional  differential equations were presented.

In this paper, we address an important topic that was not yet treated in the FC literature, namely, derivation of the closed form solutions to the Cauchy problems for the fractional differential equations with the Riemann-Liouville GFD of arbitrary order. In the 2nd section, we formulate and prove the 2nd fundamental theorem of FC for the Riemann-Liouville GFDs of arbitrary order on a suitable space of functions. This theorem provides us with the natural form of the initial conditions for the fractional differential equations with the Riemann-Liouville GFDs of arbitrary order. Then we introduce and study the sequential Riemann-Liouville GFD of arbitrary order. In the 3rd section,  we develop an operational calculus of the  Mikusi\'nski type for the GFDs of arbitrary order in the Riemann-Liouville sense. In the framework of these calculus, the GFD of the Riemann-Liouville type is reduced to a multiplication with a certain element from the field of convolution quotients. In the 4th section, this operation calculus  is applied for derivation of an explicit form of solutions to the Cauchy problems for the single- and multi-term linear fractional differential equations  with the GFDs of the Riemann-Liouville type. The solutions are obtained in form of the convolution series generated by the Sonine kernels associated to the kernels of the GFDs.

\section{GFD of arbitrary order in the Riemann-Liouville sense and some of its properties}
\label{sec2}

The GFI \eqref{GFI} and the GFD \eqref{FDR-L} of the Riemann-Liouville type with the Sonine kernels $\kappa$ and $k$ are a far reaching generalisation of the Riemann-Liouville fractional integral \eqref{RLI} and the  Riemann-Liouville fractional derivative \eqref{RLD1} that correspond to the case of the power law kernels  $\kappa(t) = h_\alpha(t)$ and $k(t)=h_{1-\alpha}(t)$ with $\alpha \in (0,1)$. Thus, the "generalized orders" of the GFI \eqref{GFI} and the GFD \eqref{FDR-L} with any Sonine kernels $\kappa$ and $k$ are restricted to the interval $(0,1)$. However, it is well-known that both the Riemann-Liouville fractional integral  and the  Riemann-Liouville fractional derivative are defined for any non-negative real order $\alpha$:
\begin{equation}
\label{RLIa}
(I^\alpha_{0+}\, f)(t) := (h_\alpha\, * \, f)(t) \, = \, \frac{1}{\Gamma(\alpha)}\, \int_0^t (t-\tau)^{\alpha-1}f(\tau)\, d\tau, \ t>0,\ \alpha >0,
\end{equation}
\begin{equation}
\label{RLDa}
(D^\alpha_{0+}\, f)(t) := \frac{d^n}{dt^n}(h_{n-\alpha}\, * \, f)(t) 
\, = \,\frac{d^n}{dt^n} (I^{n-\alpha}_{0+}\, f)(t),\ t>0,\ n\in \N,\ n-1\le \alpha < n.
\end{equation}

To define the GFI and the GFD of an arbitrary non-negative real order, in \cite{Luc21b}, the Sonine condition \eqref{Son} was extended  to the following form:
\begin{equation}
\label{Luc}
(\kappa \, * \, k)(t) = \{ 1\}^{<n>}(t),\ n\in \N,\ t>0,
\end{equation}
where
$$
\{ 1\}^{<n>}(t):= (\underbrace{\{ 1\}*\ldots\ * \{ 1\}}_{\mbox{$n$ times}})(t) = h_{n}(t) = \frac{t^{n-1}}{(n-1)!}.
$$
However, merely condition \eqref{Luc} does not ensure an essential requirement for the kernel $k$ of the GFD, namely, the condition that this kernel should have a singularity at the origin, see \cite{DGGS,Han20,HL19}. In \cite{Luc21b}, an important class of the kernels that satisfy the condition \eqref{Luc} and possess an integrable singularity of a power law type at the origin was introduced.

\begin{Definition}[\cite{Luc21b}]
\label{d_c}
Let the functions $\kappa$ and $k$ satisfy the condition \eqref{Luc} and the inclusions $\kappa \in C_{-1}(0,+\infty)$ and $k\in C_{-1,0}(0,+\infty)$ hold~true, where
\begin{equation}
\label{C-1}
C_{-1}(0,+\infty)\, := \, \{f:\ f(t)=t^{p}f_1(t),\ t>0,\ p > -1,\ f_1\in C[0,+\infty)\},
\end{equation}
\begin{equation}
\label{C-10}
 C_{-1,0}(0,+\infty) \, = \, \{f:\ f(t) = t^{p}f_1(t),\ t>0,\ -1 < p < 0,\ f_1\in C[0,+\infty)\}.
\end{equation}

The set of pairs $(\kappa,\, k)$ of such kernels is denoted by $\mathcal{L}_n$. 
\end{Definition}

In \cite{Luc21b}, the GFI and the GFD of arbitrary order with the kernels $(\kappa,\, k)\in \mathcal{L}_n$ were defined  for the first time.

\begin{Definition}[\cite{Luc21b}]
\label{dao}
Let  $(\kappa,\ k)$ be a pair of the kernels from $\mathcal{L}_n$. The~GFI with the kernel $\kappa$ is defined by the formula
\begin{equation}
\label{GFIa}
(\I_{(\kappa)}\, f)(t) :=  (\kappa\, *\, f)(t) = \int_0^t \kappa(t-\tau)f(\tau)\, 
d\tau,\ t>0.
\end{equation}
The GFD of the Riemann-Liouville type is defined as a composition of the $n$th order derivative and the GFI with the kernel $k$:
\begin{equation}
\label{FDR-La} 
(\D_{(k)}\, f)(t) := \frac{d^n}{dt^n}\,(k\, *\, f)(t) = \frac{d^n}{dt^n}\, (\I_{(k)}\, f)(t),\ t>0.
\end{equation}
Finally, the GFD of Caputo type is a regularized form of the GFD in the  Riemann-Liouville sense:
\begin{equation}
\label{FDCa}
( _*\D_{(k)}\, f)(t) := \left(\D_{(k)}\,\left(  f(\cdot) - \sum_{j=0}^{n-1}f^{(j)}(0)h_{j+1}(\cdot)\right)\right)(t),\ t>0. 
\end{equation}
\end{Definition}
In this paper, we mainly deal with the GFD of arbitrary order of the Riemann-Liouville type; the properties of the GFD of arbitrary order in the Caputo sense were investigated in \cite{Luc21b}. 

\begin{Remark}
\label{r0}
It is worth mentioning that the inclusion $(\kappa,\, k)\in \mathcal{L}_n$ can be interpreted as a statement that the GFI \eqref{GFIa} with the kernel $\kappa$ and the GFDs \eqref{FDR-La} and \eqref{FDCa} with the kernel $k$ have a "generalized order" from the interval $(n-1,\ n)$. For example, for $\alpha \in (n-1,\, n),\ n\in \N$, the kernels $\kappa = h_\alpha$ and $k = h_{n-\alpha}$ satisfy the condition \eqref{Luc}:
$$
(h_\alpha\, * \, h_{n-\alpha})(t) = h_n(t) = \{ 1\}^{<n>}(t).
$$
Moreover, the inclusions $h_\alpha \in C_{-1}(0,+\infty)$ and $h_{n-\alpha}\in C_{-1,0}(0,+\infty)$ evidently hold true. Thus,  $(h_\alpha,\, h_{n-\alpha}) \in \mathcal{L}_n$ and the GFI \eqref{GFIa} and the GFD \eqref{FDR-La} are reduced to the Riemann-Liouville fractional integral \eqref{RLIa} and the Riemann-Liouville fractional derivative \eqref{RLDa} of arbitrary non-negative order, respectively. 
\end{Remark}

In \cite{Luc21b}, the following properties of the GFI  \eqref{GFIa} of an arbitrary order on the space  $C_{-1}(0,+\infty)$ were shown:
\begin{equation}
\label{GFI-map_1}
\I_{(\kappa)}:\, C_{-1}(0,+\infty)\, \rightarrow C_{-1}(0,+\infty)\ \mbox{(mapping property)},
\end{equation}
\begin{equation}
\label{GFI-com_1}
\I_{(\kappa_1)}\, \I_{(\kappa_2)} = \I_{(\kappa_2)}\, \I_{(\kappa_1)} \ 
\mbox{(commutativity law)}, 
\end{equation}
\begin{equation}
\label{GFI-index_1}
\I_{(\kappa_1)}\, \I_{(\kappa_2)} = \I_{(\kappa_1*\kappa_2)}\ \mbox{(index law)}.
\end{equation}

As to the GFD \eqref{FDR-La} of arbitrary order, it is a left inverse operator to the GFI  \eqref{GFIa}. 

\begin{Theorem}[\cite{Luc21b}]
\label{t3_n}
Let the inclusion $(\kappa,\ k) \in \mathcal{L}_n$ hold true. 

Then the~GFD \eqref{FDR-La} is a left-inverse operator to the GFI \eqref{GFIa}  on the space $C_{-1}(0,+\infty)$:
\begin{equation}
\label{FTLn}
(\D_{(k)}\, \I_{(\kappa)}\, f) (t) = f(t),\ f\in C_{-1}(0,+\infty),\ t>0.
\end{equation}
\end{Theorem}

We mention here that the statement formulated in Theorem \ref{t3_n} is usually referred to as the 1st Fundamental Theorem of FC for the GFD \eqref{FDR-La} of arbitrary order (see \cite{Luc20} for a discussion of the 1st and the 2nd Fundamental Theorems of FC for several different kinds of the fractional derivatives). 

Moreover, it is easy to see that the GFD \eqref{FDR-La} of arbitrary order is also a right inverse operator to the GFI  \eqref{GFIa} on the space  
\begin{equation}
\label{C1k}
C_{-1,(\kappa)}^1(0,+\infty) := \{f :\,  f(t)=(\I_{(\kappa)}\, \phi)(t),\ \phi\in C_{-1}(0,+\infty)\},
\end{equation}
i.e., the relation
\begin{equation}
\label{2FTLn}
(\I_{(\kappa)}\, \D_{(k)}\, f) (t) = f(t),\ t>0
\end{equation}
holds true. 
Indeed, a function $f\in C_{-1,(\kappa)}^1(0,+\infty)$ can be represented in the form $f(t) = (\I_{(\kappa)}\, \phi)(t),\ \phi\in C_{-1}(0,+\infty)$ and thus we have the following chain of equations:
$$
(\I_{(\kappa)}\, \D_{(k)}\, f) (t) = (\I_{(\kappa)}\, \frac{d^n}{dt^n}\, (k\, *\, f) (t) = (\I_{(\kappa)}\, \frac{d^n}{dt^n}\, (k\, *\, (\kappa \, *\, \phi)) (t)=
$$
$$
(\I_{(\kappa)}\, \frac{d^n}{dt^n}\, (\{1\}^{<n>} *\, \phi)) (t)= (\I_{(\kappa)}\, \phi) (t)= f(t).
$$
It is worth mentioning that the formula \eqref{2FTLn} implicates zero initial conditions for any fractional differential equations with the GFD \eqref{FDR-La} defined on the space $C_{-1,(\kappa)}^1(0,+\infty)$ because its null space consists just of a function that is identically equal to zero for $t>0$.  

In this paper, we consider the  GFD \eqref{FDR-La} of arbitrary order on its natural domain 
\begin{equation}
\label{C1(k)}
C_{-1,(k)}^{(1)}(0,+\infty)= \{ f\in C_{-1}(0,+\infty):\, \D_{(k)}\, f \in C_{-1}(0,+\infty) \}.
\end{equation}
Evidently, the inclusion $C_{-1,(\kappa)}^{1}(0,+\infty) \subset C_{-1,(k)}^{(1)}(0,+\infty)$  holds true due to Theorem \ref{t3_n}, .
In the case of the $n$th order derivative, the space
$C_{-1,(k)}^{(1)}(0,+\infty)$ corresponds to the space of the $n$-times continuously differentiable functions, whereas the space $C_{-1,(\kappa)}^{1}(0,+\infty)$ consists of all functions that can be represented as $n$-fold integrals of continuous functions. It is worth mentioning that the spaces $C_{-1,(\kappa)}^{0}(0,+\infty)$ and $C_{-1,(k)}^{(0)}(0,+\infty)$ can be  interpreted as the space $C_{-1}(0,\, +\infty)$.

In what follows, we also employ another subspace of the space $C_{-1}(0,+\infty)$ (see \cite{Luc21c} or \cite{LucGor99} for its properties):
\begin{equation}
\label{Cn}
C_{-1}^n(0,+\infty) = \{ f\in C_{-1}(0,+\infty):\,  f^{(n)} \in C_{-1}(0,+\infty) \}.
\end{equation}

To formulate the natural initial conditions for the fractional differential equations with the GFDs of arbitrary order in the Riemann-Liouville sense defined on the space $C_{-1,(k)}^{(1)}(0,+\infty)$, we need a description of its null space that is provided by the next theorem:

\begin{Theorem}
\label{t-null}
Let  $(\kappa,\ k)$ be a pair of the kernels from $\mathcal{L}_n$ and the inclusion  $\kappa \in C_{-1}^{n-1}(0,+\infty)$ hold valid. 

Then the GFD \eqref{FDR-La} defined on the space $C_{-1,(k)}^{(1)}(0,+\infty)$ has the following $n$-dimensional null space:
\begin{equation}
\label{null}
    \mbox{Ker}\, \D_{(k)} = \left\{ \sum_{i=0}^{n-1} a_i \frac{d^{n-1-i}\kappa}{dt^{n-1-i}},\  a_{0},\dots,a_{n-1} \in \R \right\}.
\end{equation}
\end{Theorem}

\begin{proof}
The equation
$$
(\D_{(k)}\, f) (t) = \frac{d^n}{dt^n} (\I_{(k)}\, f) (t) = 0,\ t>0
$$
implicates 
$$
(\I_{(k)}\, f) (t) = P_{n-1}(t) =  a_{n-1}\frac{t^{n-1}}{(n-1)!} + a_{n-2}\frac{t^{n-2}}{(n-2)!} +\dots +
    a_{0},\ t>0,\ a_{0},\dots,a_{n-1} \in \R 
$$
that can be rewritten in the form
\begin{equation}
\label{null_1}
(k\, * \, f) (t) =  a_{n-1}\{1\}^{<n>}(t) + a_{n-2}\{1\}^{<n-1>}(t) +\dots +
    a_{0}\{1\},\ t>0.
\end{equation}
Let us now consider the convolution of the equation \eqref{null_1} with the  kernel $\kappa$ associated to the kernel $k$ of $\D_{(k)}$. Because of the condition \eqref{Luc}, the left-hand side of the equation takes the form of the $n$-fold definite integral:
$$
(\kappa \, *\,  (k\, * \, f)) (t) = ((\kappa \, * \, k)\, * \, f) (t) = (\{1\}^{<n>}\, * \, f)(t),\ t>0.
$$
The right-hand side of the equation can be represented as follows: 
$$
(\kappa \, *\, (a_{n-1}\{1\}^{<n>} + a_{n-2}\{1\}^{<n-1>} +\dots +
    a_{0}\{1\}))(t) = 
    $$
    $$
    (a_{n-1}\kappa \, * \, \{1\}^{<n>})(t) + (a_{n-2}\kappa\, *\, \{1\}^{<n-1>})(t) +\dots +
    (a_{0}\kappa\, *\, \{1\})(t).
$$
Then we arrive at the equation
$$
(\{1\}^{<n>}\, * \, f)(t) = (\{1\}^{<n>} \, *\, a_{n-1}\kappa)(t) + (\{1\}^{<n-1>} \, *\, a_{n-2}\kappa)(t) +\dots +
 (\{1\} \, *\,  a_{0}\kappa )(t).
$$
Because of the inclusion $\kappa \in C_{-1}^{n-1}(0,+\infty)$, for $t>0$, we can differentiate the last equation $n$ times and thus get the representation
$$
f(t) = a_{n-1}\kappa + a_{n-2}\frac{d\kappa}{dt} +\dots +
    a_{0}\frac{d^{n-1}\kappa}{dt^{n-1}},\ t>0 
$$
that completes the proof of theorem. 
\end{proof}

In the case of the Riemann-Liouville fractional derivative \eqref{RLDa} of  order $\alpha,\ n-1 < \alpha < n,\ n\in \N$, the kernel function $k$ is the power law function $h_{n-\alpha}$ and its associated kernel $\kappa$ is the function $h_\alpha$. Due to the evident relation $\frac{d}{dt}h_\alpha = h_{\alpha -1},\ t>0$, the formula \eqref{null} from Theorem \ref{t-null} takes the well-known form (see e.g., \cite{Samko}):
\begin{equation}
\label{null-RL}
    \mbox{Ker}\, D_{0+}^\alpha =\left\{ \sum_{i=0}^{n-1} a_i \frac{t^{\alpha-n+i}}{\Gamma(\alpha-n+i-1)},\ a_{0},\dots,a_{n-1} \in \R \right\}.
\end{equation}

The description of the null space of the GFD \eqref{FDR-La} of arbitrary order formulated in Theorem \ref{t-null} is now employed in the proof of the next important theorem:

\begin{Theorem}[2nd Fundamental Theorem of FC for the GFD of arbitrary order in the Riemann-Liouville sense]
\label{t-2ndFT}
Let  $(\kappa,\ k)$ be a pair of the kernels from $\mathcal{L}_n$ and  $\kappa \in C_{-1}^{n-1}(0,+\infty)$. 

For a function $f\in C_{-1,(k)}^{(1)}(0,+\infty)$, the formula
\begin{equation}
\label{2ndFT}
(\I_{(\kappa)}\, \D_{(k)}\, f) (t) =  f(t) - \sum_{i=0}^{n-1} \left(
\frac{d^{i}}{dt^{i}}\, \I_{(k)}\, f\right)(0)\, \frac{d^{n-1-i} \kappa}{dt^{n-1-i}},\ t>0
\end{equation}
holds valid.
\end{Theorem}

\begin{proof}
First we define an auxiliary function as follows:
\begin{equation}
\label{psi}
\psi(t) := (\I_{(\kappa)}\, \D_{(k)}\, f) (t),\ t>0.
\end{equation}
For a function $f\in C_{-1,(k)}^{(1)}(0,+\infty)$, the inclusion $\psi \in C_{-1,(\kappa)}^{1}(0,+\infty) \subset C_{-1,(k)}^{(1)}(0,+\infty)$ holds true. Thus, the GFD \eqref{FDR-La} of the function $\psi$ exists and belongs to the space $C_{-1}(0,+\infty)$. Theorem \ref{t3_n} implicates the relation
$$
(\D_{(k)}\, \psi)(t) = (\D_{(k)}\, \I_{(\kappa)}\, \D_{(k)}\, f) (t) = (\D_{(k)}\, f) (t),\ t>0
$$
that means that the function 
$\psi-f$ belongs to the kernel of the GFD $\D_{(k)}$.  According to Theorem \ref{t-null}, we get the representation
\begin{equation}
\label{F}
\psi(t)-f(t) = \sum_{i=0}^{n-1}a_{i}\frac{d^{n-1-i} \kappa}{dt^{n-1-i}},\  t>0
\end{equation}
with some constants $a_{0},\dots,a_{n-1} \in \R$. To determine these constants, we apply the GFI $\I_{(k)}$ to the representation  \eqref{F} and use the condition \eqref{Luc} for the kernels $(\kappa,\ k)\in \mathcal{L}_n$ and the inclusion  $\kappa \in C_{-1}^{n-1}(0,+\infty)$ to get the following chain of equations:
$$
(\I_{(k)}\, (\psi - f))(t)  = (\I_{(k)}\, \psi)(t) - (\I_{(k)}\, f)(t) = \left(k\, *\, 
\sum_{i=0}^{n-1}a_{i}\frac{d^{n-1-i} \kappa}{dt^{n-1-i}}\right)(t) = 
$$
$$
\sum_{i=0}^{n-1}a_{i}\, \frac{d^{n-1-i}}{dt^{n-1-i}}(k\, *\, \kappa)(t) = 
\sum_{i=0}^{n-1}a_{i}\, \{ 1\}^{<i+1>}(t) = \sum_{i=0}^{n-1}a_{i}\, \frac{t^{i}}{i!}.
$$
Combining this formula with the representation
$$
(\I_{(k)}\, \psi)(t) = (\I_{(k)}\, \I_{(\kappa)}\, \D_{(k)}\, f) (t) = 
\left( \{ 1\}^{<n>}\, * \, \D_{(k)}\, f\right) (t),
$$
we arrive at the expression
\begin{equation}
\label{psi-f}
\sum_{i=0}^{n-1}a_{i}\, \frac{t^{i}}{i!} = \left( \{ 1\}^{<n>}\, * \, \D_{(k)}\, f\right) (t) - (\I_{(k)}\, f)(t).
\end{equation}
Because of the inclusion $(\D_{(k)}\, f)\in C_{-1}(0,+\infty)$, the function 
$g(t)= (\{ 1\}^{<n>}\, * \, \D_{(k)}\, f)(t)$ belongs to the space $C^{n-1}([0,+\infty))$ and $g(0)=g^\prime(0)=\dots = g^{(n-1)}(0) = 0$ (see \cite{LucGor99}). 
Substituting $t=0$ into \eqref{psi-f}, we get a formula for the coefficient $a_0$:
\begin{equation}
\label{a0}
a_0 =  - (\I_{(k)}\, f)(0).
\end{equation}
Differentiating \eqref{psi-f} with respect to the variable $t$ leads to the relation
\begin{equation}
\label{psi-f1}
\sum_{i=1}^{n-1}a_{i}\, \frac{t^{i-1}}{(i-1)!} = \left( \{ 1\}^{<n-1>}\, * \, \D_{(k)}\, f\right) (t) - \frac{d}{dt}(\I_{(k)}\, f)(t).
\end{equation}
Substituting $t=0$ into \eqref{psi-f1}, we get a formula for the coefficient $a_1$:
\begin{equation}
\label{a1}
a_1 =  - \frac{d}{dt}(\I_{(k)}\, f)(0).
\end{equation}
Repeating the same operation $(n-1)$ times, we determine all coefficients $a_0,\, a_1,\, \dots,$ $a_{n-1}$:
\begin{equation}
\label{aj}
a_i =  - \frac{d^i}{dt^i}(\I_{(k)}\, f)(0),\ i=0,1,\dots,n-1.
\end{equation}
The formula \eqref{2ndFT}  immediately follows from the formulas \eqref{F} and \eqref{aj}  and using the representation \eqref{psi} of the auxiliary function $\psi$. 
\end{proof}

\begin{Remark}
\label{r01}
Theorem \ref{t-2ndFT} can be reformulated in terms of the projector operator of the GFD of the Riemann-Liouville type as follows:

On the space $C_{-1,(k)}^{(1)}(0,+\infty)$, the projector operator $F$ of the GFD of the Riemann-Liouville type takes the form
\begin{equation}
\label{proj}
(F\, f)(t) := f(t) - (\I_{(\kappa)}\, \D_{(k)}\, f) (t) =   \sum_{i=0}^{n-1} \left(
\frac{d^{i}}{dt^{i}}\, \I_{(k)}\, f\right)(0)\, \frac{d^{n-1-i} \kappa}{dt^{n-1-i}},\ t>0.
\end{equation}
The right-hand side of the formula \eqref{proj} specifies the natural initial conditions that should be posed while dealing with the initial-value problems for the fractional differential equations with the GFDs of arbitrary order in the Riemann-Liouville sense. Some of these problems for the linear single- and multi-term fractional differential equations will be formulated and solved in Section \ref{sec4}. 

It is worth mentioning that the function $F\, f$ always belongs to the kernel of the GFD $\D_{(k)}$:
$$
(\D_{(k)}\, (F\, f))(t) = (\D_{(k)}\, f)(t) - (\D_{(k)}\, \I_{(\kappa)}\, \D_{(k)}\, f) (t) = (\D_{(k)}\, f)(t) - (\D_{(k)}\, f)(t) = 0,\ t>0.
$$
\end{Remark}

\begin{Remark}
\label{r1}
For the GFD \eqref{FDR-La} of the "generalized order" from the interval $(0,\, 1)$ (the case of the kernels $(\kappa,\ k)\in \mathcal{L}_1$), Theorem \ref{t-2ndFT} has been formulated and proved in \cite{Luc22}. In this case, the relation  \eqref{2ndFT} takes the form 
\begin{equation}
\label{2ndFT-1}
(\I_{(\kappa)}\, \D_{(k)}\, f) (t) =  f(t) -  ( \I_{(k)}\, f)(0)\, \kappa(t),\ t>0.
\end{equation}
\end{Remark}

\begin{Remark}
\label{r2}
As mentioned in Remark \ref{r0}, the pair of the kernels  $\kappa(t)= h_{\alpha}(x),\ n-1 <\alpha <n,\ n\in \N$ and $k(t)= h_{n-\alpha}(x)$ belongs to the set of kernels $\mathcal{L}_n$. In this case, the GFI \eqref{GFIa} and the GFD \eqref{FDR-La} are reduced to the Riemann-Liouville fractional integral $I_{0+}^\alpha$  and the Riemann-Liouville fractional derivative $D_{0+}^\alpha$ of the order $\alpha$, respectively, and the formula \eqref{2ndFT} takes the well-known form (see e.g., \cite{LucSri95} or \cite{Samko}):
\begin{equation}
\label{2ndFT-RL}
(I_{0+}^\alpha\, D_{0+}^\alpha\, f) (t) = f(t) - \sum_{i=0}^{n-1} \left(
\frac{d^{i}}{dt^{i}}\, I_{0+}^{n-\alpha}\, f\right)(0)\, h_{\alpha-n+i+1}(t) = 
\end{equation}
$$
f(t) - \sum_{i=1}^{n} (D_{0+}^{\alpha-i}\, f)(0)\frac{t^{\alpha -i}}{\Gamma(\alpha -i+1)},\ t>0.
$$
\end{Remark}

In FC, one of the important research topics is the sequential fractional derivatives and the differential equations with these derivatives. The $m$-fold GFI and the  $m$-fold sequential GFDs of both the Riemann-Liouville and the Caputo types with the kernels $(\kappa,\, k) \in \mathcal{L}_1$  have been introduced and studied in \cite{Luc21d} and \cite{Luc22}.  In the rest of this section, we consider these operators in the general case of the kernels  $(\kappa,\, k) \in \mathcal{L}_n, n\in \N$. We start by defining the convolution powers of the functions and the $m$-fold GFI and the $m$-fold sequential GFD of the Riemann-Liouville type. 

\begin{Definition}
\label{d1}
Let $(\kappa,\, k) \in \mathcal{L}_n,\ n\in \N$.  The $m$-fold GFI  is defined as a composition of $m$ GFIs  with the kernel $\kappa$
\begin{equation}
\label{GFIn}
(\I_{(\kappa)}^{<m>}\, f)(t) := (\underbrace{\I_{(\kappa)} \ldots \I_{(\kappa)}}_{m\ \mbox{times}})\, f)(t) = (\kappa^{<m>}\, *\, f)(t),\ t>0,
\end{equation}
where the convolution power $g^{<m>},\ m\in \N_0$ of a function $g$ is given by the expression
\begin{equation}
\label{I-6}
g^{<m>}(t):= \begin{cases}
1,& m=0,\\
g(t), & m=1,\\
(\underbrace{g* \ldots * g}_{m\ \mbox{times}})(t),& m=2,3,\dots .
\end{cases}
\end{equation}
The $m$-fold sequential GFD in the Riemann-Liouville  sense is defined as follows:
\begin{equation}
\label{GFDLn}
(\D_{(k)}^{<m>}\, f)(t) := (\underbrace{\D_{(k)} \ldots \D_{(k)}}_{m\ \mbox{times}}\, f)(t),\  t>0.
\end{equation}
\end{Definition}

For $m=0$, we define the operators  $\I_{(\kappa)}^{<0>}$ and $\D_{(k)}^{<0>}$ as the identity operator $\mbox{Id}$.

In the formula \eqref{GFIn}, the kernel $\kappa^{<m>},\ m\in \N$  belongs to  the space
$C_{-1}(0,+\infty)$ due to Theorem \ref{t1} and thus the $m$-fold GFI  is reduced to a GFI with the kernel $\kappa^{<m>}$:
\begin{equation}
\label{GFIn-1}
(\I_{(\kappa)}^{<m>}\, f)(t)  = (\kappa^{<m>}\, *\, f)(t) = (\I_{(\kappa)^{<m>}}\, f)(t) ,\ t>0.
\end{equation}

The $m$-fold sequential GFD \eqref{GFDLn} is a generalization of the Riemann-Liouville  sequential fractional derivative to the case of the integro-differential operators with the kernels from $\mathcal{L}_{n}$.

Repeatedly applying Theorem \ref{t3_n} for the GFI \eqref{GFIa}  and the GFD \eqref{FDR-La}, we get  the following result:

\begin{Theorem}[1st Fundamental Theorem of FC for the $m$-fold sequential GFD of arbitrary order]
\label{t3_n_n}
Let $(\kappa,\, k) \in \mathcal{L}_n,\ n\in \N$. 

Then  the $m$-fold sequential GFD \eqref{GFDLn} in the Riemann-Liouville sense  is a left inverse operator to the $m$-fold  GFI \eqref{GFIn} on the space $C_{-1}(0,+\infty)$: 
\begin{equation}
\label{FTLn_n}
(\D_{(k)}^{<m>}\, \I_{(\kappa)}^{<m>}\, f) (t) = f(t),\ f\in C_{-1}(0,+\infty),\ t>0.
\end{equation}
\end{Theorem}

On a certain subspace of $C_{-1}(0,+\infty)$, the $m$-fold sequential GFD \eqref{GFDLn} of arbitrary order is also a right inverse operator to the $m$-fold GFI  \eqref{GFIn}. 

\begin{Theorem}
\label{t3_n_n1}
Let $(\kappa,\, k) \in \mathcal{L}_n,\ n\in \N$. 

Then  the $m$-fold sequential GFD \eqref{GFDLn} in the Riemann-Liouville sense  is a right inverse operator to the $m$-fold  GFI \eqref{GFIn} on the space 
\begin{equation}
\label{C1kn}
C_{-1,(\kappa)}^m(0,+\infty) := \{f :\,  f(t)=(\I_{(\kappa)}^{<m>}\, \phi)(t),\ \phi\in C_{-1}(0,+\infty)\},
\end{equation}
i.e., the relation 
\begin{equation}
\label{FTLn_n_1}
(\I_{(\kappa)}^{<m>}\, \D_{(k)}^{<m>}\,  f) (t) = f(t),\ f\in C_{-1,(\kappa)}^m(0,+\infty),\ t>0
\end{equation}
holds true. Moreover, on this space, the $m$-fold sequential GFD \eqref{GFDLn} can be represented as a GFD with the kernel $k^{<m>}$:
\begin{equation}
\label{GFD_a_n}
(\D_{(k)}^{<m>}\,  f) (t) = (\D_{(k)^{<m>}}\,  f) (t),\ f\in C_{-1,(\kappa)}^m(0,+\infty),\ t>0.
\end{equation}
\end{Theorem}

\begin{proof}
For $f\in C_{-1,(\kappa)}^m(0,+\infty)$, we have  a representation $f(t) = (\I_{(\kappa)}^{<m>}\, \phi)(t),$ $\phi\in C_{-1}(0,+\infty)$ that leads to the following chain of equations due to Theorem \ref{t3_n_n}:
$$
(\I_{(\kappa)}^{<m>}\, \D_{(k)}^{<m>}\, f) (t) = (\I_{(\kappa)}^{<m>}\, \D_{(k)}^{<m>} \I_{(\kappa)}^{<m>}\, \phi)(t)=
(\I_{(\kappa)}^{<m>}\, \phi)(t) = f(t).
$$
This proves the relation \eqref{FTLn_n_1}. As to the formula \eqref{GFD_a_n}, we just evaluate its left- and right-hand sides for a function from $C_{-1,(\kappa)}^m(0,+\infty)$ using Theorems \ref{t3_n} and \ref{t3_n_n} and the representation \eqref{GFIn-1}:
$$
(\D_{(k)}^{<m>}\, f) (t) =( \D_{(k)}^{<m>}\, \I_{(\kappa)}^{<m>}\, \phi)(t) = \phi(t),\ t>0
$$
and
$$
(\D_{(k)^{<m>}}\, f) (t) =( \D_{(k)^{<m>}}\, \I_{(\kappa)}^{<m>}\, \phi)(t) =
( \D_{(k)^{<m>}}\, \I_{(\kappa)^{<m>}}\, \phi)(t) = \phi(t),\ t>0
$$
that completes the proof of the theorem. 
\end{proof}

Of course, the formula \eqref{FTLn_n_1} does not hold true for the $m$-fold sequential GFD \eqref{GFDLn}  defined on its natural domain 
\begin{equation}
\label{C1(k)_n}
C_{-1,(k)}^{(m)}(0,+\infty)= \{ f\in C_{-1}(0,+\infty):\, \D_{(k)}^{<i>}\, f \in C_{-1}(0,+\infty),\ i=1,\dots,m \}.
\end{equation}
For a pair of the kernels $(\kappa,\, k) \in \mathcal{L}_n,\ n\in \N$, the space $C_{-1,(\kappa)}^{m}(0,+\infty)$ defined by \eqref{C1kn} is a subspace of the space $ C_{-1,(k)}^{(m)}(0,+\infty)$ because of Theorem \ref{t3_n_n}. On the space $C_{-1,(k)}^{(m)}(0,+\infty)$, the projector operator of the $m$-fold sequential GFD \eqref{GFDLn} takes the form described in the next theorem. 

\begin{Theorem}
\label{t-proj_n}
Let $f \in C_{-1,(k)}^{(m)}(0,+\infty)$.  The projector operator $F_m$ of the $m$-fold sequential GFD \eqref{GFDLn} can be represented as follows:
\begin{equation}
\label{proj-n}
(F_m\, f)(t) := f(t) - (\I_{(\kappa)}^{<m>}\, \D_{(k)}^{{<m>}}\, f) (t) = \sum_{i=0}^{m-1} (\I_{(\kappa)}^{<i>}\, (F\, (\D_{(k)}^{<i>)}\, f))(t),\ t>0,
\end{equation}
where $F$ stands for the projector operator \eqref{proj} of the GFD $\D_{(k)}$ of the Riemann-Liouville type. 
\end{Theorem}

\begin{proof}
On the space $C_{-1,(k)}^{(m)}(0,+\infty)$, the identity operator $\mbox{Id}$ can be represented as follows:
$$
\mbox{Id} = (\mbox{Id}-\I_{(\kappa)}\, \D_{(k)})+(\I_{(\kappa)}\, \D_{(k)}-\I_{(\kappa)}^{<2>}\, \D_{(k)}^{{<2>}})+\dots
+(\I_{(\kappa)}^{<m-1>}\, \D_{(k)}^{{<m-1>}}-\I_{(\kappa)}^{<m>}\, \D_{(k)}^{{<m>}})+
$$
$$
+\I_{(\kappa)}^{<m>}\, \D_{(k)}^{{<m>}} =
(\mbox{Id}-\I_{(\kappa)}\, \D_{(k)})+(\I_{(\kappa)}\, (\D_{(k)}-(\I_{(\kappa)}\, \D_{(k)})\, \D_{(k)}))+\dots
$$
$$
+(\I_{(\kappa)}^{<m-1>}\, (\D_{(k)}^{{<m-1>}}-(\I_{(\kappa)}\, \D_{(k)})\, \D_{(k)}^{{<m-1>}}))
+\I_{(\kappa)}^{<m>}\, \D_{(k)}^{{<m>}} \, = \, F+(\I_{(\kappa)}\, F(\D_{(k)})) +\dots
$$
$$
+(\I_{(\kappa)}^{<m-1>}\, F(\D_{(k)}^{{<m-1>}}))+\I_{(\kappa)}^{<m>}\, \D_{(k)}^{{<m>}}.
$$
The formula \eqref{proj-n} immediately follows from the last representation. 
\end{proof}

Evidently, the function $(F_m\, f)(t)$ belongs to the kernel of the  $m$-fold sequential GFD \eqref{GFDLn}:
$$
(\D_{(k)}^{{<m>}}\, (F_m\, f))(t) = (\D_{(k)}^{{<m>}}\, f)(t) - ((\D_{(k)}^{{<m>}}\, \I_{(\kappa)}^{<m>}\, \D_{(k)}^{{<m>}}\, f) (t) = 
$$
$$
(\D_{(k)}^{{<m>}}\, f)(t) - (\D_{(k)}^{{<m>}}\, f)(t) = 0.
$$
On the other hand, the inclusions
$$
\mbox{Ker}\, \D_{(k)} \subset  \mbox{Ker}\, \D_{(k)}^{{<2>}} \subset \dots \subset \mbox{Ker}\, \D_{(k)}^{{<m>}}
$$
follow directly from the definition of the $m$-fold sequential GFD \eqref{GFDLn}. 

In the rest of the paper, we deal with the single- and multi-term fractional differential equations with the $m$-fold sequential GFDs \eqref{GFDLn} in the Riemann-Liouville sense. The formulas \eqref{proj}  and \eqref{proj-n} for the projector operators of the GFD $\D_{(k)}$ of the Riemann-Liouville type and the $m$-fold sequential GFD $\D_{(k)}^{{<m>}}$ allow us to formulate the natural initial conditions for these fractional differential equations. The solution method is based on the Mikusi\'nski type operational calculus for the  GFD $\D_{(k)}$ of the Riemann-Liouville type that will be developed in the next section. 

\section{Operational calculus for the GFD of arbitrary order in the Riemann-Liouville sense}
\label{sec3}

In the recent paper \cite{Luc21c}, a Mikusi\'nski type operational calculus  for the GFD  in the Caputo sense with the "generalized order" from the interval $(0,\, 1)$ (the case of the kernels $(\kappa,\ k)\in \mathcal{L}_1$) was constructed. In this section, we discuss a modification and an extension of this operational calculus to the case  of the GFD   of arbitrary order in the Riemann-Liouville sense  (the case of the kernels $(\kappa,\ k)\in \mathcal{L}_n,\ n\in \N$). 

In what follows, we deal with the GFI defined by \eqref{GFIa} and the GFD of arbitrary order in the Riemann-Liouville sense defined by \eqref{FDR-La}. The kernels $\kappa$ and $k$  of the GFI and the GFD, respectively, are from the set $\mathcal{L}_n$ and thus they  satisfy the condition \eqref{Luc}. 

Because the GFI  has the same form for the GFDs both in the Riemann-Liouville sense and in the Caputo sense, some results presented in \cite{Luc21c} remain valid for the operational calculus for the GFD of arbitrary order in the Riemann-Liouville sense. In particular, we mention the following important theorem:

\begin{Theorem}[\cite{LucGor99}]
\label{t1}
The   triple   $\mathcal{R}_{-1} = (C_{-1}(0,+\infty),+,*)$     with  the  usual
addition $+$ and  multiplication $*$ in form of  the  Laplace convolution is a commutative ring without divisors of zero.
\end{Theorem}

Thus, the~GFI with the kernel $\kappa$ can be interpreted as a multiplication on the ring $\mathcal{R}_{-1}$:
\begin{equation}
\label{GFIar}
(\I_{(\kappa)}\, f)(t) =  (\kappa\, *\, f)(t),\ t>0.
\end{equation}

Because the  GFD of arbitrary order in the Riemann-Liouville sense is a left inverse operator to the GFI \eqref{GFIa}  (Theorem \ref{t3_n})  and  there exists no unity element with respect to multiplication in  $\mathcal{R}_{-1}$, representing GFD as a  multiplication on the ring $\mathcal{R}_{-1}$ is not possible (see \cite{Luc21c} for details). However, following a standard procedure, we extend the ring $\mathcal{R}_{-1}$ to a field of convolution quotients and define an algebraic GFD as a multiplication with a special element of this field.

The convolution quotients are equivalence classes on the set $C_{-1}^2(0,+\infty):=$ $ C_{-1}(0,+\infty) \times  (C_{-1}(0,+\infty) \setminus \{0\})$ 
that are generated by the natural equivalence
relation in the form
$$
(f_1,\, g_1)\sim     (f_2,\, g_2)    \Leftrightarrow     (f_1\, * \,
g_2)(t) = (f_2\, * \,  g_1)(t),\ (f_1,\, g_1),\ (f_2,\, g_2) \in C_{-1}^2(0,+\infty).
$$
The equivalence classes $C_{-1}^2(0,+\infty)/\sim$ are usually denoted as quotients:
$$
\frac{f}{g}:= \{ (f_1,\, g_1) \in C_{-1}^2(0,+\infty):\ (f_1,\, g_1) \sim (f,\, g) \}.
$$
On the set of equivalence classes, the usual operations of addition and multiplication are introduced:
$$
\frac{f_1}{g_1} + \frac{f_2}{g_2}:= \frac{ f_1\, *\,  g_2\, +\, f_2\, * g_1}
{g_1\, *\, g_2},
$$
$$
\frac{f_1}{g_1}\cdot  \frac{f_2}{g_2} :=\frac{f_1\, *\,  f_2}{g_1\, *\,  g_2}.
$$
It is easy to verify that the operations $+$ and $\cdot$ on $C_{-1}^2(0,+\infty)/\sim$ are correctly defined because they do not depend on the representatives of the  equivalence classes.

The following important Theorem is an immediate implication of Theorem \ref{t1}:

\begin{Theorem}[\cite{LucGor99}]
\label{t8}
The triple $\mathcal{F}_{-1} = (C_{-1}^2(0,+\infty)/\sim,\ +,\ \cdot)$ is a  field that is usually referred to as the field of convolution quotients.
\end{Theorem}

It is worth mentioning that the ring $\mathcal{R}_{-1}$ is embedded into the field
$\mathcal{F}_{-1}$:
\begin{equation}
\label{emb}
f  \mapsto \frac{f\, *\, \kappa}{\kappa},
\end{equation}
where $\kappa$ is the kernel of the GFI \eqref{GFIa}. Of course, in the formula \eqref{emb}, $\kappa$ can be replaced by any other non-zero element of the space $C_{-1}(0,+\infty)$.

Another useful operation on $C_{-1}^2(0,+\infty)/\sim$ is multiplication with a scalar $\lambda \in \R$ or $\lambda \in \Com$:
$$
\lambda\, \frac{f}{g}:=\frac{\lambda\, f}{g},\ \frac{f}{g}\in C_{-1}^2(0,+\infty)/\sim.
$$

Evidently, the space of functions $C_{-1}(0,+\infty)$ is a vector space that induces a vector space structure on the set  $C_{-1}^2(0,+\infty)/\sim$ of equivalence classes  with the operations "+" and multiplication with a scalar. Because  the constant function $\{ \lambda \}$ 
is an element of the ring  $\mathcal{R}_{-1}$,  one has to
distinguish between multiplication with a scalar $\lambda$ on the vector
space  $C_{-1}^2(0,+\infty)/\sim$ and  multiplication with the constant
function $\{ \lambda \}$ in the field $\mathcal{F}_{-1}$ of convolution quotients that is defined as follows:
$$
\{\lambda\}\cdot \frac{f}{g} = \frac{\{\lambda\}\, *\,  f}{g},\ \frac{f}{g}\in \mathcal{F}_{-1}.
$$

As mentioned in \cite{Luc21c}, the unity element $I = \frac{\kappa}{\kappa}$  of the field $\mathcal{F}_{-1}$ with respect to multiplication  does not belong to the ring $\mathcal{R}_{-1}$. Thus, it is not a conventional function but rather a kind of a generalized function (hyperfunction in the terminology of \cite{Yos})  that in our operational calculus plays the role of the Dirac $\delta$-function.  

According to the formula \eqref{GFIar}, the GFI \eqref{GFIa} can be interpreted as multiplication with the element $\kappa\in \mathcal{R}_{-1}$. The embedding \eqref{emb} assigns to $\kappa$ the  element $\frac{\kappa\, *\, \kappa}{\kappa} \in \mathcal{F}_{-1}$. Its inverse element is another important hyperfunction that will be used in this section for defining an algebraic GFD of the Riemann-Liouville type. 

\begin{Definition}[\cite{Luc21c}]
\label{d4}
The convolution quotient
\begin{equation}
\label{alg}
S_\kappa = \frac{\kappa}{\kappa\, *\, \kappa} \in \mathcal{F}_{-1}
\end{equation}
is called an algebraic inverse element to the GFI \eqref{GFIa}.
\end{Definition}

Evidently, the relation 
\begin{equation}
\label{alg-inv}
{\kappa}\cdot S_\kappa  = I
\end{equation}
holds true, where $I$ is the unity of $\mathcal{F}_{-1}$ with respect to multiplication. It justifies the next definition.

\begin{Definition}
\label{d5}
Let $(\kappa,\, k) \in \mathcal{L}_n,\ n\in \N$, the inclusion $\kappa \in C_{-1}^{n-1}(0,\, +\infty)$ hold true, and  $f\in C_{-1}(0,\, +\infty)$. The algebraic GFD of arbitrary order in the Riemann-Liouville sense  is defined as follows:
\begin{equation}
\label{AGFD}
\D_{(k)}\, f = S_\kappa \cdot f - S_\kappa \cdot (F\, f),
\end{equation}
where the function $f$ and the projector operator $F\, f$ defined by the formula \eqref{proj} are interpreted as elements of the field $\mathcal{F}_{-1}$.
\end{Definition}

In the previous section, the GFD \eqref{FDR-La} of arbitrary order in the Riemann-Liouville sense was defined on the space of functions $C_{-1,(k)}^{(1)}(0,+\infty)$ that consists of all functions from $C_{-1}(0,+\infty)$ whose GFD exists and belongs to the space $C_{-1}(0,+\infty)$. As we see from Definition \ref{d5}, the algebraic GFD \eqref{AGFD} is a generalized derivative that assigns a certain element from the field  $\mathcal{F}_{-1}$ to any function $f$ from $C_{-1}(0,+\infty)$.  Thus, the algebraic GFD \eqref{AGFD} makes sense also for the functions whose GFD \eqref{FDR-La} does not exist in the usual sense. However, for any function  $f\in C_{-1,(k)}^{(1)}(0,+\infty)$, the GFD \eqref{FDR-La} of arbitrary order in the Riemann-Liouville sense coincides with the algebraic GFD \eqref{AGFD}. That's why we use the same notations for both derivatives even if their domains are very different. 

\begin{Theorem}
\label{t9}
Let  $(\kappa,\, k)\in  \mathcal{L}_n$ and the inclusion $\kappa \in C_{-1}^{n-1}(0,\, +\infty)$ hold true. For a function $f\in C_{-1,(k)}^{(1)}(0,+\infty)$, its algebraic GFD \eqref{AGFD} coincides with the GFD \eqref{FDR-La} of arbitrary order in the Riemann-Liouville sense:
\begin{equation}
\label{AGFD-f}
(\D_{(k)}\, f)(t) = \D_{(k)}\, f  = S_\kappa \cdot f - S_\kappa \cdot (F\, f),
\end{equation}
where the projector operator $F\, f$ is defined by the formula \eqref{proj} and the functions $f$, $\D_{(k)}\, f$, and $F\, f$ are interpreted as elements of the field $\mathcal{F}_{-1}$.
\end{Theorem}

\begin{proof}
For a function $f\in C_{-1,(k)}^{(1)}(0,+\infty)$, the formula \eqref{proj} deduced in the previous section leads to the following representation:
\begin{equation}
\label{AGFD-f-1}
(\I_{(\kappa)}\, \D_{(k)}\, f) (t)  =  (\kappa \, *\, (\D_{(k)}\, f))(t) = f(t) - (F\, f)(t), \ t>0. 
\end{equation}
Evidently, the right- and the left-hand sides of the relation \eqref{AGFD-f-1} are from the space $C_{-1,(k)}^{(1)}(0,+\infty) \subset C_{-1}(0,\, +\infty)$ and thus we can interpret it as an equality of two elements from the convolution quotients field $\mathcal{F}_{-1}$. Now we multiply this equality with the element $S_\kappa$ and use the relation \eqref{alg-inv} to arrive at the representation \eqref{AGFD-f}.
\end{proof}

\begin{Remark}
\label{r-RL}
A formula of type \eqref{AGFD-f} for the Riemann-Liouville fractional derivative of arbitrary non-negative order $\alpha,\ n-1\le \alpha < n,\, n\in \N$ has been derived in \cite{Luc93} for the first time. As mentioned in Remark \ref{r2}, in this formula, the projector operator \eqref{proj} for the Riemann-Liouville fractional derivative has the well-known form
$$
(F\, f)(t) = \sum_{i=1}^{n} (D_{0+}^{\alpha-i}\, f)(0)\frac{t^{\alpha -i}}{\Gamma(\alpha -i+1)},\ t>0.
$$
Thus, the natural initial condition for the fractional differential equations with the  Riemann-Liouville fractional derivative $D_{0+}^\alpha$ is as follows:
\begin{equation}
\label{R-L-ini}
(F\, f)(t) =  \sum_{i=1}^{n} (D_{0+}^{\alpha-i}\, f)(0)\frac{t^{\alpha -i}}{\Gamma(\alpha -i+1)} = \gamma (t),\ \gamma(t) \in \mbox{Ker} (D_{0+}^\alpha).
\end{equation}
Because of the representation 
$$
\mbox{Ker} (D_{0+}^\alpha) =\left\{ \sum_{i=1}^{n} a_i\frac{t^{\alpha -i}}{\Gamma(\alpha -i+1)},\ a_1,\dots,a_n \in \R\right\},
$$
the initial condition \eqref{R-L-ini} can be rewritten in the standard form:
\begin{equation}
\label{R-L-ini-st}
 (D_{0+}^{\alpha-i}\, f)(0) = a_i,\ i=1,\dots,n.
\end{equation}
\end{Remark}

Now we proceed with an extension of Definition \ref{d5}  to the case of the $m$-fold sequential GFD \eqref{GFDLn}.

\begin{Definition}
\label{d6}
Let $(\kappa,\, k) \in \mathcal{L}_n,\ n\in \N$, the inclusion $\kappa \in C_{-1}^{n-1}(0,\, +\infty)$ hold true, and  $f\in C_{-1}(0,\, +\infty)$. The $m$-fold sequential algebraic GFD of arbitrary order in the Riemann-Liouville sense  is defined as follows:
\begin{equation}
\label{AGFDn}
\D_{(k)}^{<m>}\, f = S_\kappa^m \cdot f - S_\kappa^m \cdot (F_m\, f),
\end{equation}
where the function $f$ and  the projector operator $F_m\, f$ defined by the formula \eqref{proj-n} are interpreted as elements of the field $\mathcal{F}_{-1}$.
\end{Definition}

The $m$-fold sequential algebraic GFD of arbitrary order in the Riemann-Liouville sense is well  defined for any function from $C_{-1}(0,\, +\infty)$. Thus, it can be interpreted as a kind of a generalized derivative. However, for the functions from the space
$C_{-1,(k)}^{(m)}(0,+\infty)$, Theorem \ref{t-proj_n} and the  arguments we employed in the case $m=1$ ensure the result formulated in the next theorem.

\begin{Theorem}
\label{t10}
Let $(\kappa,\, k) \in \mathcal{L}_n,\ n\in \N$ and the inclusion $\kappa \in C_{-1}^{n-1}(0,\, +\infty)$ hold true. For a function $f\in C_{-1,(k)}^{(m)}(0,+\infty)$, its $m$-fold sequential algebraic GFD \eqref{AGFDn} coincides with the $m$-fold sequential GFD \eqref{GFDLn}:
\begin{equation}
\label{AGFDn-f}
(\D_{(k)}^{<m>}\, f)(t) = \D_{(k)}^{<m>}\, f  = S_\kappa^m \cdot f - S_\kappa^m \cdot  (F_m\, f),
\end{equation}
where the projector operator $F_m$ is defined by the formula \eqref{proj-n} and and the functions $f$, $\D_{(k)}^{<m>}\, f$, and $F_m\, f$ are interpreted as elements of the field $\mathcal{F}_{-1}$
\end{Theorem}

\begin{Remark}
\label{rproj}
Taking into account the formula \eqref{proj-n} for the projector operator $F_m$ and the relations $(I_{(\kappa)}^{<i>}\, f)(t) = (\kappa^{<i>}\, *\, f)(t)$ and $S_\kappa \cdot \kappa = I$, the representation \eqref{AGFDn-f} can be rewritten in a form more convenient for dealing with the fractional differential equations containing the GFDs:
\begin{equation}
\label{form1}
(\D_{(k)}^{<m>}\, f)(t) = \D_{(k)}^{<m>}\, f  = S_\kappa^m \cdot f - \sum_{i=0}^{m-1}S_\kappa^{m-i} \cdot (F\,(\D_{(k)}^{<i>} f)),
\end{equation}
where the projector operator $F$ is given by the expression \eqref{proj}. 
\end{Remark}

The formulas \eqref{AGFD-f} and \eqref{AGFDn-f} (or \eqref{form1}) mean that  on the field $\mathcal{F}_{-1}$ of convolution quotients, the GFD  \eqref{FDR-La}  and the $m$-fold sequential GFD \eqref{GFDLn} are reduced to multiplication  with certain field elements. Moreover, the initial conditions in terms of the corresponding  projector operators are also integrated into these formulas. Using the formulas \eqref{AGFD-f} and \eqref{AGFDn-f} or \eqref{form1} we can rewrite the initial-value problems for the linear fractional differential equations with  the GFDs  \eqref{FDR-La}  and the $m$-fold sequential GFDs \eqref{GFDLn}  as some algebraic (in fact, linear) equations on the field $\mathcal{F}_{-1}$. The solutions to these equations can be easily obtained in explicit form.  However, in general, they  are elements of the field $\mathcal{F}_{-1}$ and thus  generalized functions or hyperfunctions. The embedding \eqref{emb} of the ring  $\mathcal{R}_{-1}$ into the field $\mathcal{F}_{-1}$ ensures that some hyperfunctions can be  interpreted as conventional functions from the space $C_{-1}(0,\, +\infty)$. In the rest of this section, we introduce an important class of such hyperfunctions in form of rational functions $R(S_\kappa) = Q(S_\kappa)/P(S_\kappa)$, where $Q$ and $P$ are polynomials and $\mbox{deg}(Q) < \mbox{deg}(P)$. We start with a result regarding the so-called convolution series generated by the functions from the space $C_{-1}(0,+\infty)$. 

\begin{Theorem}[\cite{Luc22}]
\label{t11}
Let a function $\kappa \in C_{-1}(0,+\infty)$  be represented in the form
\begin{equation}
\label{rep}
\kappa(t) = h_{p}(t)\kappa_1(t),\ t>0,\ p>0,\ \kappa_1\in C[0,+\infty)
\end{equation} 
and the convergence radius of the power  series
\begin{equation}
\label{ser}
\Sigma(z) = \sum^{+\infty }_{j=0}a_{j}\, z^j,\ a_{j}\in \Com,\ z\in \Com
\end{equation}
be non-zero. Then the convolution series 
\begin{equation}
\label{conser}
\Sigma_\kappa(t) = \sum^{+\infty }_{j=0}a_{j}\, \kappa^{<j+1>}(t)
\end{equation}
is convergent for all $t>0$ and defines a function  from the space $C_{-1}(0,+\infty)$. 
Moreover, the series 
\begin{equation}
\label{conser_p}
t^{1-\alpha}\, \Sigma_\kappa(t) = \sum^{+\infty }_{j=0}a_{j}\, t^{1-\alpha}\, \kappa^{<j+1>}(t),\ \ \alpha = \min\{p,\, 1\}
\end{equation}
is uniformly convergent for $t\in [0,\, T]$ for any $T>0$.
\end{Theorem}

The convolution series \eqref{conser} is a far reaching generalization of the power series \eqref{ser} that corresponds to the case of the kernel $\kappa = \{1\}$ (Mikusi\'nski's operational calculus for the first order derivative).

\begin{Example}
\label{ex1}
For the geometric series
\begin{equation}
\label{geom}
\Sigma(t) = \sum_{j=1}^{+\infty} \lambda^{j-1}t^j,\ \lambda \in \Com,\ \lambda\not = 0,\ t\in \Com
\end{equation}
with the convergence radius $r = 1/|\lambda|$, the convolution series \eqref{conser} takes the form
\begin{equation}
\label{l}
l_{\kappa,\lambda}(t) = \sum_{j=1}^{+\infty} \lambda^{j-1}\kappa^{<j>}(t),\ \lambda \in \Com,\ t>0.
\end{equation}
According to Theorem \ref{t11}, this series is convergent for all $t>0$ and defines a function from the space $C_{-1}(0,+\infty)$.
\end{Example} 

\begin{Remark}
The  Mikusi\'nski operational calculus for the first order derivative (\cite{Mik59}) is based on the kernel function $\kappa = \{ 1\}$. Because $\kappa^{<j>}(t) = \{ 1\}^{<j>}(t) = h_{j}(t)$, the convolution series \eqref{l} takes the form of the exponential function:
\begin{equation}
\label{l-Mic}
l_{\kappa,\lambda}(t) = \sum_{j=1}^{+\infty} \lambda^{j-1}h_j(t) =
\sum_{j=0}^{+\infty} \frac{(\lambda\, t)^j}{j!} = e^{\lambda\, t}.
\end{equation}
\end{Remark}

\begin{Remark}
For constructing the operational calculus of Mikusi\'nski type for the Riemann-Liouville fractional derivative (\cite{LucSri95}, \cite{HadLuc}),  the kernel $\kappa = h_{\alpha}$ of the Riemann-Liouville fractional integral has been used. Due to the relation $\kappa^{<j>}(t) = h_{\alpha}^{<j>}(t) = h_{j\alpha}(t)$, the convolution series \eqref{l} has the form
\begin{equation}
\label{l-Cap}
l_{\kappa,\lambda}(t) = \sum_{j=1}^{+\infty} \lambda^{j-1}h_{j\alpha}(t) =
t^{\alpha-1}\sum_{j=0}^{+\infty} \frac{\lambda^j\, t^{j\alpha}}{\Gamma(j\alpha+\alpha)} = t^{\alpha -1}E_{\alpha,\alpha}(\lambda\, t^{\alpha}),
\end{equation}
where the two-parameters Mittag-Leffler function $E_{\alpha,\beta}$ is defined by the following absolutely convergent series:
\begin{equation}
\label{ML}
E_{\alpha,\beta}(z) = \sum_{j=0}^{+\infty}\frac{z^j}{\Gamma(\alpha\, j + \beta)},\ \Re(\alpha)>0,\ z,\beta \in \Com.
\end{equation}
\end{Remark}

For other important particular cases of the convolution series \eqref{l} we refer to the recent publication \cite{Luc21c}. 

The role of the function $l_{\kappa,\lambda}$ in the operational calculus for the GFD of arbitrary order in the Riemann-Liouville sense is illustrated in the next theorem. 

\begin{Theorem}[\cite{Luc21c}]
\label{t12}
For any  element of the convolution quotients field $\mathcal{F}_{-1}$ in form $(S_\kappa -\lambda) $ with $\lambda \in \Com$,   its inverse element with respect to multiplication is provided by the  convolution series  $l_{\kappa,\lambda} \in \mathcal{R}_{-1}\subset \mathcal{F}_{-1}$  defined by \eqref{l} i.e., the operational relation 
\begin{equation}
\label{op-rel}
\frac{I}{S_{\kappa} - \lambda}\, = \, l_{\kappa,\lambda}(t),\ t>0
\end{equation}
holds true.
\end{Theorem}

In \cite{Luc21c}, a proof of this theorem was presented for the case of the kernels $(\kappa,\, k)\in \mathcal{L}_{1}$. However, in the general case of the kernels $(\kappa,\, k)\in \mathcal{L}_{n},\ n\in \N$, all arguments from this proof and thus the statement of Theorem \ref{t12} remain valid. 

\begin{Example}
\label{ex2}
For the kernel  $\kappa =  \{1\}$  (Mikusi\'nski's operational calculus for the first order derivative), the operational relation \eqref{op-rel} takes the form (see \eqref{l-Mic})
\begin{equation}
\label{l-Mic-op}
\frac{I}{S_{\kappa} - \lambda}\, = \, l_{\kappa,\lambda}(t) \, = \,  e^{\lambda\, t}.
\end{equation}
\end{Example}

\begin{Example}
\label{ex3}
For the kernel   $\kappa(t) = h_{\alpha}(t),\ t>0$ (operational calculi for the Riemann-Liouville and  the Caputo fractional derivatives \cite{HadLuc,Luc93,LucGor99,LucSri95}), the operational relation \eqref{op-rel} has the form (see \eqref{l-Cap})
\begin{equation}
\label{l-Cap-op}
\frac{I}{S_{\kappa} - \lambda}\, = \, l_{\kappa,\lambda}(t) \, = \, t^{\alpha -1}E_{\alpha,\alpha}(\lambda\, t^{\alpha}),
\end{equation}
where the two-parameters Mittag-Leffler function $E_{\alpha,\beta}$ is defined by \eqref{ML}.
\end{Example}

For other important particular cases of the operational relation \eqref{op-rel} see \cite{Luc21c}. 

Based on the formula \eqref{op-rel}, other useful operational relations can be easily deduced in terms of the convolution powers of the function  $l_{\kappa,\lambda}$ defined by the formula \eqref{l}. Because convolution of the functions from the ring $\mathcal{R}_{-1}$ complies with  multiplication of the corresponding elements of the field $\mathcal{F}_{-1}$, we immediately get the following operational relation:
\begin{equation}
\label{op-rel-m}
\frac{I}{(S_{\kappa} - \lambda)^m}\, = \, l^{<m>}_{\kappa,\lambda}(t),\ t>0,\ m\in \N.
\end{equation}
It is worth mentioning that the convolution powers $l^{<m>}_{\kappa,\lambda}$ can be represented as some convolution series by means of  the Cauchy product for the convolution series \eqref{l} (see \cite{Luc21c} for more details). 

\begin{Example}
\label{ex4}
In the case of the kernel $\kappa = \{1\}$ (Mikusi\'nski's operational calculus for the first order derivative), we get the well-known operational relation (see \cite{Mik59}):
\begin{equation}
\label{l-Mic-op-m}
\frac{I}{(S_{\kappa} - \lambda)^m}\, = \,  h_m(t)\, e^{\lambda\, t}.
\end{equation}
\end{Example}

\begin{Example}
\label{ex5}
For the kernel $\kappa(t) = h_{\alpha}(t),\ t>0$ (operational calculi for  the Riemann-Liouville and the Caputo fractional derivatives), the operational relation \eqref{op-rel-m} takes the form (\cite{LucGor99}):
\begin{equation}
\label{l-Cap-op-m}
\frac{I}{(S_{\kappa} - \lambda)^m}\, = \, t^{m\alpha -1}E_{\alpha,m\alpha}^m(\lambda t^\alpha),\ t>0,\ m\in \N,
\end{equation}
where the Mittag-Leffler type function $E_{\alpha,\beta}^{m}$ is defined as follows:
$$
E_{\alpha,\beta}^{m}(z):=
\sum^{\infty
}_{j=0}\frac{(m)_j z^{j}}{  j!
\Gamma (\alpha j + \beta)}, \ \alpha,\beta>0, \ z\in \Com, \
(m)_j = \prod_{i=0}^{j-1} (m+i).
$$
\end{Example}

Finally, we mention that the operational relation \eqref{op-rel-m} allows us to represent any rational function in $S_{\kappa}$ with the denominator's degree grater than the numerator's degree as a conventional function from the space $C_{-1}(0,+\infty)$. 

\begin{Theorem}
\label{trational}
Let $R(S_\kappa) = Q(S_\kappa)/P(S_\kappa)$, where $Q$ and $P$ are polynomials and $\mbox{deg}(Q) < \mbox{deg}(P)$. Then the operational relation
\begin{equation}
\label{rational}
R(S_\kappa) = \sum_{j=1}^{J}\sum_{i=1}^{m_j} a_{ij}\, l^{<i>}_{\kappa,\lambda_j}(t),\ t>0,\ \sum_{j=1}^J m_j = \mbox{deg}\, (P),
\end{equation}
holds true, where the constants $\lambda_j$ and $m_j$, $j=1,\dots,J$ are uniquely determined by representation of the rational function $R(S_\kappa)$ as a sum of the partial fractions:
\begin{equation}
\label{R-partial}
R(S_\kappa) = \sum_{j=1}^{J}\sum_{i=1}^{m_j} \frac{a_{ij}}{(S_\kappa - \lambda_j)^i},\ \sum_{j=1}^J m_j = \mbox{deg}\, (P).
\end{equation}
\end{Theorem}

The proof of Theorem \ref{t13} immediately follows from the representation \eqref{R-partial} and the operational relation \eqref{op-rel-m}.

\section{Fractional differential equations with the GFDs of arbitrary order in the Riemann-Liouville sense}
\label{sec4}

Whereas the fractional differential equations, both ordinary and partial, with the GFDs in the Caputo sense have been discussed in several publications 
(see e.g., \cite{Koch11,LucYam16,Koch19_1,KK,LucYam20,Sin18,Sin20,Luc21c}), to the best of author's knowledge, only one publication dealing with the fractional differential equations with the GFDs in the Riemann-Liouville sense has been published until now. In \cite{Luc21d}, analytical formulas for solutions to some linear single- and multi-term ordinary differential equations with the GFDs in the Riemann-Liouville sense of the generalized order from the interval $(0,\, 1)$ i.e., with the kernels from the set $\mathcal{L}_1$ have been derived by using the method of convolution series. In this section, we consider the initial-value problems for the differential equations with the GFDs of arbitrary order in the Riemann-Liouville sense with the suitably formulated initial conditions. Their solutions will be derived by using the Mikusi\'nski type operational calculus developed in the previous section.

To demonstrate the method, we start with a simple equation that however has not yet been considered in the literature. Let $(\kappa,\, k) \in \mathcal{L}_n,\ n\in \N$, the inclusion $\kappa \in C_{-1}^{n-1}(0,\, +\infty)$ hold true, and  $f\in C_{-1}(0,\, +\infty)$. We consider the following Cauchy problem:
\begin{equation}
\label{eq-1-1}
\begin{cases}
(\D_{(k)}\, y)(t) - \lambda y(t) = f(t), & \lambda \in \R,\ t>0, \\
(F\, y)(t) = \gamma_0(t), & \gamma_0(t)\in \mbox{Ker}(\D_{(k)}),
\end{cases}
\end{equation}
where the projector operator $F$ is defined by \eqref{proj} and the null space of $\D_{(k)}$ is given by the formula \eqref{null}.

In \eqref{eq-1-1}, $\D_{(k)}$ stands for the  GFD of arbitrary order in the Riemann-Liouville sense and the unknown function $y$ is looked for in the space $C_{-1,(k)}^{(1)}(0,+\infty)$. 

Comparing the formulas \eqref{null} and \eqref{proj}, we can represent the initial condition from \eqref{eq-1-1} in a more usual form:
\begin{equation}
\label{eq-1-1-f}
\begin{cases}
(\D_{(k)}\, y)(t) - \lambda y(t) = f(t), & \lambda \in \R,\ t>0, \\
\left(\frac{d^{i}}{dt^{i}}\, \I_{(k)}\, f\right)(0) = a_i,& i=0,\dots,n-1.
\end{cases}
\end{equation}

However, the expression $(F\, y)(t) = \gamma_0(t), \ \gamma_0(t)\in \mbox{Ker}(\D_{(k)})$ provides a natural form of the initial condition that can be generalized to the case of the fractional differential equations with the $m$-fold sequential GFDs of arbitrary order in the Riemann-Liouville sense. 

For a function $y\in C_{-1,(k)}^{(1)}(0,+\infty)$, the  GFD of arbitrary order in the Riemann-Liouville sense coincides with the algebraic GFD on the field $\mathcal{F}_{-1}$ of convolution quotients (see  the formula \eqref{AGFD-f} from Theorem \ref{t9}). Thus,  the initial-value problem \eqref{eq-1-1} is reduced to the following linear equation on $\mathcal{F}_{-1}$:
\begin{equation}
\label{eq-1-2}
S_\kappa \cdot y - S_\kappa \cdot \gamma_0  - \lambda y = f,
\end{equation}
where the functions $y,\, f$ and $\gamma_0$  from the ring  $\mathcal{R}_{-1}$ are interpreted as the elements of the field  $\mathcal{F}_{-1}$.

The unique solution to the linear equation \eqref{eq-1-2} takes the form
\begin{equation}
\label{eq-1-3}
y  = f \cdot \frac{I}{S_\kappa - \lambda} + \gamma_0  \cdot \frac{S_\kappa}{S_\kappa - \lambda}.
\end{equation}

The next step is an interpretation of the solution \eqref{eq-1-3}  as a conventional function. First we represent it in the form
\begin{equation}
\label{eq-1-4}
y(t)  = y_f(t) + y_{iv}(t),\ y_f:= f \cdot \frac{I}{S_\kappa - \lambda},\ y_{iv}:= \gamma_0  \cdot \frac{S_\kappa}{S_\kappa - \lambda}.
\end{equation}

The operational relation \eqref{op-rel} from Theorem \ref{t12} and the embedding of the ring $\mathcal{R}_{-1}$ into the field $\mathcal{F}_{-1}$ immediately provide us with a suitable interpretation of the term $y_f$:
\begin{equation}
\label{eq-1-5}
y_f (t) = f \cdot \frac{I}{S_\kappa - \lambda} = (f\, *\, l_{\kappa,\lambda})(t) =
\int_0^t l_{\kappa,\lambda}(\tau)f(t-\tau)\, d\tau,
\end{equation}
where the function $l_{\kappa,\lambda}$ is the convolution series \eqref{l}. 

The part $y_{iv}$ of the solution also can be interpreted as a conventional function from the ring $\mathcal{R}_{-1}$: 
\begin{equation}
\label{op-rel-1}
y_{iv}= \gamma_0  \cdot \frac{S_\kappa}{S_\kappa - \lambda}  = \gamma_0  \cdot \frac{S_\kappa -\lambda +\lambda }{S_\kappa - \lambda} = \gamma_0  \cdot \left( I + \frac{\lambda }{S_\kappa - \lambda}\right) = \gamma_0 + \lambda (\gamma_0 \, *\, l_{\kappa,\lambda})(t).
\end{equation}

Now we come back to the representation $\gamma_0 (t) = \sum_{i=0}^{n-1}
    a_{i}\frac{d^{i}\kappa}{dt^{i}}$ of the function $\gamma_0 \in \mbox{Ker}(\D_{(k)})$ from the initial condition in \eqref{eq-1-1}. The formula \eqref{proj} for the projector operator $F$ leads then to the initial conditions in form of $n$ equations $\left(\frac{d^{i}}{dt^{i}}\, \I_{(k)}\, f\right)(0) = a_i, i=0,\dots,n-1$ (see the formula \eqref{eq-1-1-f}). 

 Making use of the identity $\kappa \cdot S_\kappa  = I$ and  the inclusion $\kappa \in C_{-1}^{n-1}(0,\, +\infty)$,  the part $y_{iv}$ of the solution can be represented as an element of the ring $\mathcal{R}_{-1}$ as follows:
\begin{equation}
\label{sol-1-y-p}
y_{iv}(t) = \gamma_0  \cdot \frac{S_\kappa}{S_\kappa - \lambda} = \sum_{i=0}^{n-1} 
    a_{i}\frac{d^{i}\kappa}{dt^{i}} \cdot \frac{S_\kappa}{S_\kappa - \lambda} = 
    \sum_{i=0}^{n-1} 
    a_{i}\frac{d^{i}}{dt^{i}} \frac{\kappa \cdot S_\kappa}{S_\kappa - \lambda} = 
    \end{equation} 
    $$
   \sum_{i=0}^{n-1} 
    a_{i}\frac{d^{i}}{dt^{i}} \frac{I}{S_\kappa - \lambda} =  \sum_{i=0}^{n-1} 
    a_{i}\frac{d^{i}}{dt^{i}}  l_{\kappa,\lambda}(t),
 $$
where $l_{\kappa,\lambda}$ is the convolution series \eqref{l}.  

Summarizing the derivations above, we get a proof of the following theorem:

\begin{Theorem}
\label{t13}
Let $(\kappa,\, k) \in \mathcal{L}_n,\ n\in \N$, the inclusion $\kappa \in C_{-1}^{n-1}(0,\, +\infty)$ hold true, and  $f\in C_{-1}(0,\, +\infty)$. The unique solution $y\in C_{-1,(k)}^{(1)}(0,+\infty)$  to the initial-value problem \eqref{eq-1-1} for the fractional differential equation with the GFD of arbitrary order in the Riemann-Liouville sense is given by the expression
\begin{equation}
\label{sol-1}
y (t) = y_f(t) + y_{iv}(t) = (f\, *\, l_{\kappa,\lambda})(t) + \gamma_0 + \lambda (\gamma_0 \, *\, l_{\kappa,\lambda})(t),
\end{equation}
where $l_{\kappa,\lambda}$ is the convolution series \eqref{l}.

For the initial conditions as in equation \eqref{eq-1-1-f}, this  solution can be represented as follows:
\begin{equation}
\label{sol-1-p}
y (t) = y_f(t) + y_{iv}(t) = (f\, *\, l_{\kappa,\lambda})(t) + 
\sum_{i=0}^{n-1} 
    a_i \, \frac{d^{i}}{dt^{i}}  l_{\kappa,\lambda}(t).
\end{equation}

The function $y_f$ is the solution to the problem with the homogeneous initial condition ($\gamma_0(t) \equiv 0$), whereas the function $y_{iv}$ is the solution to the problem with the homogeneous equation ($f(t) \equiv 0$). 
\end{Theorem}

\begin{Remark}
\label{rem2}
The Cauchy problem
\eqref{eq-1-1} with the GFD of the generalized order from the interval $(0,\, 1)$ (the case of the kernels $(\kappa,\, k) \in \mathcal{L}_1$) has been solved in \cite{Luc21d} by employing the method of convolution series. In this case, 
the function  $\gamma_0 \in \mbox{Ker}(\D_{(k)})$ has the form
$\gamma_0 (t) =   a_0\kappa(t)$. According to the formula \eqref{sol-1-y-p}, the part $y_{iv}$ of the solution can be represented as follows:
\begin{equation}
\label{sol-1-y-p-1}
y_{iv}(t) = a_{0} l_{\kappa,\lambda}(t).
\end{equation}  
Thus, the solution \eqref{sol-1-p} to the Cauchy problem
\eqref{eq-1-1-f} with the initial condition in the form $\left( \I_{(k)}\, y\right)(0) = a_0$  is given by the formula
$$
y (t) =  (f\, *\, l_{\kappa,\lambda})(t) + a_{0} l_{\kappa,\lambda}(t),
$$
that is exactly the same formula that was derived in \cite{Luc21d}.
\end{Remark}

\begin{Example}
\label{ex-RL-1}
In the case of the Riemann-Liouville fractional derivative of order $\alpha,\ n-1<\alpha < n,\ n\in \N$ ($\kappa(t) = h_{\alpha}(t),\ k(t) = h_{n-\alpha}(t)$), the function $l_{\kappa,\lambda}$ is given by the formula \eqref{l-Cap} in terms of the two-parameters Mittag-Leffler function. Because this power series is absolutely and uniformly convergent on any closed interval $[\epsilon,\, T]$, we can differentiate it term by term and get the formula
\begin{equation}
\label{l-Cap-dif}
\frac{d^{i}}{dt^{i}}  l_{\kappa,\lambda}(t) = \frac{d^{i}}{dt^{i}}  t^{\alpha -1}E_{\alpha,\alpha}(\lambda\, t^{\alpha}) =  t^{\alpha -1-i}E_{\alpha,\alpha-i}(\lambda\, t^{\alpha}),\ t>0.
\end{equation}
Theorem \ref{t13} and the representation \eqref{sol-1-y-p}  lead to the well-known solution formula
\begin{equation}
\label{sol-1-2}
y (t) = (f(\tau)\, *\, \tau^{\alpha -1}E_{\alpha,\alpha}(\lambda\, \tau^{\alpha}))(t) + \sum_{i=0}^{n-1} a_i\, t^{\alpha -1-i}E_{\alpha,\alpha-i}(\lambda\, t^{\alpha})
\end{equation}
to the following Cauchy problem for the one-term fractional differential equation with the Riemann-Liouville fractional derivative of order $\alpha$ with $\alpha \in (n-1,\, n),\ n\in \N$:
\begin{equation}
\label{eq-1-1-f-RL}
\begin{cases}
(\D_{0+}^\alpha\, y)(t) - \lambda y(t) = f(t), & \lambda \in \R,\ t>0, \\
\left(\D_{0+}^{\alpha-n+i} \, f\right)(0) = a_i,& i=0,\dots,n-1.
\end{cases}
\end{equation}
In \cite{LucSri95}, this formula was deduced by means of an operational calculus for the Riemann-Liouville fractional derivative. 
\end{Example}

In the rest of this section, we consider a class of the multi-term fractional differential equations with the sequential GFDs of arbitrary order in the Riemann-Liouville sense. 

Let $P_m(z) =\sum^{m}_{p=0}c_{p}z^{p}$ be a polynomial with the real or 
  complex
coefficients,  $(\kappa,\, k) \in \mathcal{L}_n,\ n\in \N$, the inclusion $\kappa \in C_{-1}^{n-1}(0,\, +\infty)$ hold true, and  $f\in C_{-1}(0,\, +\infty)$. In analogy to the one-term case, we formulate the following Cauchy problem for a linear multi-term fractional differential equation with the sequential GFDs and the constant coefficients:
\begin{equation}
\label{(20.1)}
\begin{cases}
\sum_{p=0}^m c_{p}(\D_{(k)}^{<p>}\, y)(t)  = f(t), & c_p \in \Com,\ p=0,1,\dots,m,\,  t>0, \\
(F\, \D_{(k)}^{<p>})(t) = \gamma_p(t), & \gamma_p(t)\in \mbox{Ker}(\D_{(k)}),\ p=0,\dots, m-1,
\end{cases}
\end{equation}
where $\D_{(k)}^{<p>}$ stands for the  sequential GFD of arbitrary order in the Riemann-Liouville sense, the projector operator $F$ is defined by the equation \eqref{proj}, the null space of $\D_{(k)}$ is given by the formula \eqref{null}, and the unknown function $y$ is looked for in the space $C_{-1,(k)}^{(m)}(0,+\infty)$. 

Our main result concerning the  problem \eqref{(20.1)} is formulated in the next theorem:

\begin{Theorem}
\label{t(20.1)}
The unique solution  to  the  Cauchy 
problem \eqref{(20.1)} for the multi-term fractional differential equation with the GFDs in the Riemann-Liouville sense 
 on the space $C_{-1,(k)}^{(m)}(0,+\infty)$ can be represented as follows:
\begin{equation}
\label{(20.5)}
y(t)= y_f(t) + y_{iv}(t), \ t>0,
\end{equation}
where $y_f$ is the unique  solution to the Cauchy 
problem \eqref{(20.1)} with the homogeneous initial conditions ($\gamma_p(t) \equiv 0,\ p=0,\dots,m-1$) in the form 
\begin{equation}
\label{(20.5-1)}
 y_f(t) = \sum^{J}_{j=1}\sum^{m_{j}}_{i=1}c_{ij}(f\, * \, l_{\kappa,\lambda_j}^{<i>})(t),\ \sum^{J}_{j=1} m_j = m,
\end{equation}
and $y_{iv}$ is the unique solution to the Cauchy 
problem \eqref{(20.1)} for the homogeneous equation  ($f(t) \equiv 0$) in the form 
\begin{equation}
\label{(20.5-2)}
y_{iv} = \gamma_{0}(t) +
\sum^{m-1}_{p=0}\sum^{J}_{j=1}\sum^{m_{j}}_{i=1}c_{ijp}(\gamma_p \, * \, l_{\kappa,\lambda_j}^{<i>})(t),\ \sum^{J}_{j=1} m_j = m.
\end{equation}
In the formulas \eqref{(20.5-1)} and \eqref{(20.5-2)}, the convolution series $l_{\kappa,\lambda}$ is defined by the equation \eqref{l} and the constants $\lambda_j$ , $c_{ij}$, $c_{ijp}$  are
determined by the representations
of the rational functions $\frac{1}{P_m(z)}$ and $\frac{P_p(z)}{P_m(z)}$, $P_p(z) =\sum^{m-p}_{j=1}c_{p+j}z^{j}$,  $p=0,\dots,m-1$ 
as  the  sums  of  the partial
fractions
\begin{equation}
\label{(20.6)}
\frac{1}{P_m(z)} = \frac{1}{c_{m}(z-\lambda_{1})^{m_1}\times \ldots
\times (z-\lambda_{J})^{m_J}}= \sum^{J}_{j=1}\sum^{m_{j}}_{i=1}\frac{c_{ij}}{(z-\lambda_{j})^{i}},
\end{equation}
\begin{equation}
\label{(20.6)-1}
\frac{P_p(z)}{P_m(z)} =
\sum^{J}_{j=1}\sum^{m_{j}}_{i=1}\frac{c_{ijp}}{(z-\lambda_{j})^{i}},
p=1,..,m-1; \ \frac{P_0(z)}{P_m(z)}=1+\sum^{J}_{j=1}\sum^{m_{j}}_{i=1}
\frac{c_{ij0}}{(z-\lambda_{j})^{i}}.
\end{equation}
\end{Theorem}

\begin{proof}

Our strategy for solving the Cauchy problem \eqref{(20.1)} is the same as the one for the one-term problem \eqref{eq-1-1}.

First, we apply the formula \eqref{form1} to rewrite the Cauchy problem \eqref{(20.1)} as an algebraic (in fact,  linear) equation on the field $\mathcal{F}_{-1}$ of convolution quotients: 
\begin{equation}
\label{form2}
c_0 y + \sum_{p=1}^m c_{p} (S_\kappa^p \cdot y - \sum_{j=0}^{p-1}S_\kappa^{p-j} \cdot \gamma_p) = f,
\end{equation}
where the functions $f,\, y,$ and $\gamma_p,\ p=0,\dots m-1$ are interpreted as elements of the field $\mathcal{F}_{-1}$. The equation \eqref{form2} can be rewritten in the form 
\begin{equation}
\label{form3}
P_m(S_\kappa) \cdot y = f + \sum^{m-1}_{p=0} \gamma_p \cdot P_p(S_\kappa), \ P_p (S_\kappa) =\sum^{m-p}_{j=1}c_{p+j}S_{\kappa}^{j},\ p=0,1,\dots,m-1. 
\end{equation}
The linear equation \eqref{form3} has a unique solution on the field $\mathcal{F}_{-1}$ of convolution quotients:
\begin{equation}
\label{(20.4)}
y = y_f + y_{iv} = f\cdot \frac{I} {P_m(S_\kappa)} +\sum^{m-1}_{p=0} \gamma_p \cdot \frac{P_p(S_\kappa)}{ P_m(S_\kappa)}, \ P_p (S_\kappa) =\sum^{m-p}_{j=1}c_{p+j}S_{\kappa}^{j}.
\end{equation}

In general, the right-hand side of the formula \eqref{(20.4)} is an element of the field $\mathcal{F}_{-1}$, i.e. a hyperfuntion. Because the  ring $\mathcal{R}_{-1}$  is  embedded  into  the  field $\mathcal{F}_{-1}$  of
convolution quotients and the linear equation \eqref{form3}  possesses a unique solution in form \eqref{(20.4)},   the Cauchy 
problem
(\ref{(20.1)}) has a unique solution on the space $C_{-1,(k)}^{(m)}(0,+\infty)$ if and only if the right-hand side of the formula 
(\ref{(20.4)})  can be interpreted as a conventional function from the space $C_{-1,(k)}^{(m)}(0,+\infty)$. Let us check that it is really the case. 

First we mention that the solution formulas \eqref{(20.5)}, \eqref{(20.5-1)}, and \eqref{(20.5-2)} immediately follow from the representation \eqref{(20.4)}, Theorem \ref{trational}, and the relations \eqref{(20.6)} and \eqref{(20.6)-1}. The only delicate case is the rational function $\frac{P_0(S_\kappa)}{P_m(S_\kappa)}$ because $\mbox{deg}(P_0) = \mbox{deg}(P_m)= m$. Thus, we have to proceed as in the case of the one-term equation \eqref{eq-1-1} and to represent this rational function as follows:
$$
\frac{P_0(S_\kappa)}{P_m(S_\kappa)} = \frac{P_0(S_\kappa) + c_0 -c_0}{P_m(S_\kappa)} = \frac{P_m(S_\kappa) -c_0}{P_m(S_\kappa)} = 
I - \frac{c_0}{P_m(S_\kappa)}.
$$
Using this representation and the formula \eqref{(20.6)-1}, we arrive at the relation
\begin{equation}
\label{p=0}
\gamma_0 \cdot \frac{P_0(S_\kappa)}{P_m(S_\kappa)}  = \gamma_0 \cdot \left(I - \frac{c_0}{P_m(S_\kappa)}\right) = \gamma_0(t) + \sum^{J}_{j=1}\sum^{m_{j}}_{i=1}c_{ij0}
(\gamma_0 \, * \, l_{\kappa,\lambda_j}^{<i>})(t)
\end{equation}
that is exactly  the part of the formula \eqref{(20.5-2)} that corresponds to the rational function $\frac{P_0(S_\kappa)}{P_m(S_\kappa)}$. 

Now we prove  that the right-hand side of the formula 
(\ref{(20.4)})  can be interpreted as a conventional function from the space $C_{-1,(k)}^{(m)}(0,+\infty)$. 

We start with mentioning that any element of $\mathcal{F}_{-1}$ in form $\frac{I} {S_\kappa-\lambda},\ \lambda \in \Com$ is reduced
to a conventional function from 
the space $C_{-1}(0,\, +\infty)$ (see Theorem \ref{t12}). By definition, this space  coincides with the space $C_{-1,(k)}^{(0)}(0,+\infty)$. 

To proceed, we prove now the  following  auxiliary statement: if
a function $u$ is from the space $C_{-1,(k)}^{(s)}(0,+\infty),\ s\in \N_0$, then
the function
\begin{equation}
\label{(20.7)}
v = u \cdot \frac{I} {S_\kappa-\lambda} =
u \, * \, l_{\kappa, \lambda}
\end{equation}
belongs to the space $C_{-1,(k)}^{(s+1)}(0,+\infty)$. Indeed,
since $u$ and $l_{\kappa, \lambda}$ are from the space  $C_{-1}(0,\, +\infty)$, the inclusion 
$v  \in   C_{-1}(0,\, +\infty)$   holds true because of Theorem \ref{t1}.   Multiplying the equation \eqref{(20.7)} with ${S_\kappa-\lambda}$ and then with the kernel $\kappa$ and using the identity $\kappa \cdot S_\kappa = I$, we arrive at the following equation for the function $v$: 
\begin{equation}
\label{(20.8)}
v(t) = (\I_{(\kappa)}\, u)(t)+\lambda (\I_{(\kappa)}\, v)(t), 
\end{equation}
where $\I_{(\kappa)}$ is the GFI with the kernel $\kappa$. 
Because the functions $u$ and $v$ are from the space $C_{-1}(0,\, +\infty)$, Theorem \ref{t3_n} leads to the inclusion $v \in C_{-1,(k)}^{(1)}(0,+\infty)$. If the function $u$ is from the space $C_{-1,(k)}^{(s)}(0,+\infty)$ with $s=0$, the auxiliary statement is proved. If $s\ge 1$, the functions $u$ and $v$ belong to the space $ C_{-1,(k)}^{(1)}(0,+\infty)$ and the  relation (\ref{(20.8)}) along with Theorem \ref{t3_n} ensures the inclusion $v \in C_{-1,(k)}^{(2)}(0,+\infty)$. Repeating the same arguments again and again, we arrive at the inclusion $v \in C_{-1,(k)}^{(s+1)}(0,+\infty)$ that we wanted to prove. 

Now we employ the inclusion $f\in C_{-1}(0,\, +\infty)$, the 
representation \eqref{(20.6)} with $m_{1}+\ldots
+m_{J}=m$, and the 
auxiliary statement proved above to  get  the
inclusion $y_f = f\cdot \frac{I} {P_m(S_\kappa)} \in  C_{-1,(k)}^{(m)}(0,+\infty)$. As to the part $y_{iv}$ of the solution, we first remind the inclusions $\gamma _{p} \in  \ker  \D_{(k) }$, $p=0,1,\dots,m-1$. Consequently, 
the functions $\gamma_{p},\ p=0,1,\dots,m-1$ belong  to the space  $C_{-1,(k)}^{(s)}(0,+\infty)$ for any
 $s\in \N$. For $p=1,\dots,m-1$, the inclusions $\gamma_p \cdot \frac{P_p(S_\kappa)}{ P_m(S_\kappa)} \in C_{-1,(k)}^{(s)}(0,+\infty)$ with
arbitrary $s\in \N$ immediately follow from the inclusions $\gamma_p \in C_{-1,(k)}^{(s)}(0,+\infty)$ and the auxiliary statement proved above. 
For $p=0$, we refer to the representation (see the formula \eqref{p=0}):
$$
\gamma_0 \cdot \frac{P_0(S_\kappa)}{P_m(S_\kappa)}  = \gamma_0 \cdot \left(I - \frac{c_0}{P_m(S_\kappa)}\right) = \gamma_0(t) - \gamma_0 \cdot \frac{c_0}{P_m(S_\kappa)}
$$
that leads to the inclusion $\gamma_0 \cdot \frac{P_0(S_\kappa)}{P_m(S_\kappa)} \in C_{-1,(k)}^{(s)}(0,+\infty)$ with
arbitrary $s\in \N$ by the same arguments that we employed above. 
\end{proof}

\begin{Remark}
The initial conditions $(F\, \D_{(k)}^{<p>})(t) = \gamma_p(t), \ \gamma_p(t)\in \mbox{Ker}(\D_{(k)}),$ $p=0,1,\dots,m-1$ from the Cauchy problem \eqref{(20.1)} can be represented in a more usual form (see the formula \eqref{proj} for the projector operator $F$):
\begin{equation}
\label{(20.1)-1}
\begin{cases}
\sum_{p=0}^m c_{p}(\D_{(k)}^{<p>}\, y)(t)  = f(t), & c_p \in \Com,\ p=0,1,\dots,m,\,  t>0, \\
\left(\frac{d^{i}}{dt^{i}}\, \I_{(k)}\, \D_{(k)}^{<p>} \right)(0) = a_{ip},& i=0,\dots,n-1,\  p=0,\dots, m-1.
\end{cases}
\end{equation}
\end{Remark}

\begin{Remark}
The solution formulas  of type \eqref{(20.5)}, \eqref{(20.5-1)}, and \eqref{(20.5-2)} to the Cauchy problem \eqref{(20.1)} for the multi-term fractional differential equations with the Riemann-Liouville fractional derivatives ($\kappa(t) = h_{\alpha}(t),\ k(t) = h_{n-\alpha}(t)$) were derived in \cite{Luc93} and  \cite{LucSri95} by using an operational calculus of the  Mikusi\'nski type for the Riemann-Liouville fractional derivative. 
\end{Remark}

\end{document}